%% file: GT1.tex
\def\selecta{n}
\def\draft{n}

\documentstyle[12pt,amssymb,epic,eepic]{amsart}		

\input macros.tex
\input defs.tex

\begin{document}

\title[Associators and the Grothendieck-Teichmuller Group]
  {On Associators and the Grothendieck-Teichmuller Group I}

\author{Dror Bar-Natan}
\address{Institute of Mathematics\\
        The Hebrew University\\
        Giv'at-Ram, Jerusalem 91904\\
        Israel}
\email{drorbn@math.huji.ac.il}

\thanks{This article is available electronically at
  {\tt http://www.ma.huji.ac.il/$\sim$drorbn}, and at \newline
  {\tt http://xxx.lanl.gov/abs/q-alg/9606021}.
}

\ifselecta{}{\dedicatory{
  Appeared in Selecta Mathematica, New Series {\bf 4} (1998) 183--212
}}
\ifselecta{}{\date{This edition: Oct.~19,~1998; \ \ First edition: June 21, 1996.}}

\begin{abstract}
We present a formalism within which the relationship (discovered
by Drinfel'd in~\cite{Drinfeld:QuasiHopf, Drinfeld:GalQQ}) between
associators (for quasi-triangular quasi-Hopf algebras) and (a variant of)
the Grothendieck-Teichmuller group becomes simple and natural, leading
to a simplification of Drinfel'd's original work. In particular, we
reprove that rational associators exist and can be constructed
iteratively, though the proof itself still depends on the apriori
knowledge that a not-necessarily-rational associator exists.
\end{abstract}

\maketitle

\ifselecta{}{\tableofcontents}

\input{intro}
\input{definitions}
\input{ASS}
\input{theorem}
\input{Just}
\input{refs}

\end{document}
\endinput

%% file: macros.tex
\theoremstyle{plain}
\newtheorem{principle}{Principle}

\newtheorem{theorem}{Theorem}

\newtheorem{proposition}{Proposition}[section]
\newtheorem{lemma}[proposition]{Lemma}
\newtheorem{corollary}[proposition]{Corollary}
\newtheorem{claim}[proposition]{Claim}

\theoremstyle{definition}
\newtheorem{definition}[proposition]{Definition}

\theoremstyle{remark}

\newtheorem{remark}[proposition]{Remark}
\newtheorem{warning}[proposition]{Warning}

\def\printname#1{
	\if\draft y
		\smash{\makebox[0pt]{\hspace{-0.5in}
			\raisebox{8pt}{\tt\tiny #1}}}
	\fi
}

\newcommand{\mathmode}[1]{$#1$}

\newlength{\standardunitlength}
\catcode`\@=11
\long\def\@makecaption#1#2{%
    \vskip 10pt
    \setbox\@tempboxa\hbox{
      \small\sf{\bfcaptionfont #1. }\ignorespaces #2}%
    \ifdim \wd\@tempboxa >\captionwidth {%
        \rightskip=\@captionmargin\leftskip=\@captionmargin
        \unhbox\@tempboxa\par}%
      \else
        \hbox to\hsize{\hfil\box\@tempboxa\hfil}%
    \fi}
\font\bfcaptionfont=cmssbx10 scaled \magstephalf
\newdimen\@captionmargin\@captionmargin=2\parindent
\newdimen\captionwidth\captionwidth=\hsize
\catcode`\@=12

\newcommand{\Aut}{\operatorname{Aut}}

\newcommand{\mor}{\operatorname{mor}}

\def\lbl#1{\label{#1}\printname{#1}}

\def\qed{{\hfill\text{$\Box$}}}
\def\eqdef{\overset{\text{def}}{=}}

\newenvironment{myitemize}{
	\begin{list}{$\bullet$}{\setlength{\leftmargin}{16pt}
        \setlength{\labelwidth}{12pt}
        \setlength{\labelsep}{4pt}}
}{
	\end{list}
}
\newenvironment{myenumi}{
	\begin{list}{\theunumi}{\usecounter{enumi}
	\setlength{\leftmargin}{20pt}
        \setlength{\labelwidth}{16pt}
        \setlength{\labelsep}{4pt}}
}{
	\end{list}
}

\newlength{\globalparindent}

%% file: defs.tex
\advance \topmargin by -\headheight
\advance \topmargin by -\headsep
\evensidemargin \oddsidemargin
\newcommand{\G}[1]{{\cal G}_{#1}}
\newcommand{\F}[1]{{\cal F}_{#1}}

\def\KZ{\text{\it KZ}}

\def\Apb{{\cal A}^{pb}}
\def\Apbc{\hat{\cal A}^{pb}}
\def\Apbm{{\cal A}^{pb(m)}}

\def\ASSc{\widehat{\text{{\bf ASS}}}}
\def\ASSm{\text{{\bf ASS}}^{(m)}}
\def\ASSmm{\text{{\bf ASS}}^{(m-1)}}
\def\GT{\text{{\bf GT}}}
\def\GTc{\widehat{\text{{\bf GT}}}}
\def\GTm{\text{{\bf GT}}^{(m)}}
\def\GRT{\text{{\bf GRT}}}
\def\GRTc{\widehat{\text{{\bf GRT}}}}
\def\GRTm{\text{{\bf GRT}}^{(m)}}
\def\GRTmm{\text{{\bf GRT}}^{(m-1)}}
\def\PB{\text{{\bf PaB}}}
\def\PBc{\widehat{\text{{\bf PaB}}}}
\def\PBm{\text{{\bf PaB}}^{(m)}}
\def\PCD{\text{{\bf PaCD}}}
\def\PCDc{\widehat{\text{{\bf PaCD}}}}
\def\PCDm{\text{{\bf PaCD}}^{(m)}}
\def\PP{\text{{\bf PaP}}}
\def\TMP{\text{{\bf TMP}}}
\def\Zc{{\widehat{\text{{\bf Z}}}}}
\def\Zm{\text{{\bf Z}}^{(m)}}

\def\pentagon{{
  \setlength{\unitlength}{6\standardunitlength}
  \begin{array}{c}  \hspace{-1.7mm}
     \raisebox{-4pt}{
       \begin{picture}(20,35)(0,-10)
         \path(5,0)(15,0)(20,10)(10,20)(0,10)(5,0)
       \end{picture}
     }
     \hspace{-1.9mm}
  \end{array}
}}

\def\hexagon{{
  \setlength{\unitlength}{6\standardunitlength}
  \begin{array}{c}  \hspace{-1.7mm}
     \raisebox{-4pt}{
       \begin{picture}(18,35)(0,-10)
         \path(9,20)(0,15)(0,5)(9,0)(18,5)(18,15)(9,20)
       \end{picture}
     }
     \hspace{-1.9mm}
  \end{array}
}}

\newcommand\ifselecta[2]{\if\selecta y #1 \else #2 \fi}

%% file: intro.tex
\section{Introduction}

\subsection{Reminders about quasi-triangular quasi-Hopf algebras}
A quasi-triangular quasi-Hopf algebra \cite{Drinfeld:QuasiHopf}
is an algebra $A$ together with a not-quite-cocommutative and
not-quite-coassociative coproduct $\Delta$, whose failure to be
cocommutative is ``controlled'' by some element $R\in A^{\otimes 2}$
and whose failure to be coassociative is ``controlled'' by some element
$\Phi\in A^{\otimes 3}$ (for more details, see~\cite{Drinfeld:QuasiHopf}
or~\cite{Kassel:Book, ShniderSternberg:Book}). For the representations of
$A$ to form a tensor category, $R$ and $\Phi$ have to obey the so-called
``pentagon'' $\pentagon\!$ and ``hexagon'' $\hexagon\!_\pm$ equations
(see section~\ref{ASS}). In~\cite{Drinfeld:QuasiHopf} Drinfel'd
finds a ``universal'' formula $(R_{\KZ},\Phi_{\KZ})$ for a solution
of $\pentagon\!$ and $\hexagon\!_\pm$ by considering holonomies of
the so-called Knizhnik-Zamolodchikov connection. The formula $R_{\KZ}$
is very simple --- $R_{\KZ}$ is in a clear sense ``an exponential''. The
formula $\Phi_{\KZ}$ is somewhat less satisfactory, as it requires
analysis --- differential equations and/or iterated integrals whose
values are most likely transcendental numbers~\cite{Drinfeld:QuasiHopf,
LeMurakami:Universal, Zagier:Values}. In~\cite{Drinfeld:GalQQ}
Drinfel'd proves that there is an iterative algebraic procedure for
finding a universal formula for a solution $(R,\Phi)$ of $\pentagon\!$,
$\hexagon\!_\pm$ (with $R=R_{\KZ}$), and that such a universal formula
(called {\em an associator}) can be found iteratively and over the
rationals.

Associators (and the iterative procedure for constructing them) are
important in the theory of finite-type invariants of knots (Vassiliev
invariants) \cite{Bar-Natan:VasBib, Bar-Natan:NAT, Cartier:Construction,
Kassel:Book, LeMurakami:Universal, Piunikhin:Combinatorial} and of
3-manifolds \cite{LeMurakamiOhtsuki:Universal, Le:UniversalIHS}.
Recently, Etingof and Kazhdan~\cite{EtingofKazhdan:Quantization,
EtingofKazhdan:Poisson} used associators to show that any Lie bialgebra
can be quantized. Their results become algorithmically computable once
we know that an associator can be found iteratively.\footnote{
  For most applications of associators to finite-type invariants, it is
  in fact sufficient to use a weaker but more complicated notion of an
  associator for which an iterative construction
  was given in~\cite{Bar-Natan:NAT}.
}

Unfortunately, Drinfel'd's paper is complicated and hard to read. It
involves the introduction, almost ``out of thin air'', of two groups,
$\GTc$ and $\GRTc$, that act on the set $\ASSc$ of all associators. Both
groups act simply transitively on $\ASSc$, with $\GTc$ acting on the right
and $\GRTc$ on the left, and the two actions commute. He then studies
these groups and their actions on $\ASSc$ to deduce the existence of
formulae better then $\Phi_{\KZ}$. Drinfel'd's ``{\bf G}rothendieck-{\bf
T}eichmuller'' group $\GTc$ is closely related to number theory and the
group $\text{Gal}(\bar{{\Bbb Q}}/{\Bbb Q})$.  See~\cite{Drinfeld:GalQQ,
Schneps:DessinsdEnfants}. $\GRTc$ is in some sense a ``g{\bf R}aded''
version of $\GTc$, explaining why Drinfel'd inserted an R in the middle
of its name.

\subsection{What we do} The purpose of this paper is to present a
framework within which the set of associators $\ASSc$, the groups $\GTc$
and $\GRTc$, and the relevant facts about them are natural.  In fact,
the mere fact that $\GTc$ and $\GRTc$ exist and act simply transitively
on the right (for $\GTc$) and on the left (for $\GRTc$), with the two
actions commuting, stems from the following basic principle (which I
learned from M.~Hutchings):

\begin{principle} \lbl{BasicPrinciple}
If $B$ is a mathematical structure (i.e., a set, a set
with a basepoint, an algebra, a category, etc.) and if $C$ is an
isomorphic mathematical structure, then on the set $A$ of all isomorphisms
$B\to C$ there are two commuting group actions, with both actions
simple and transitive:
\begin{itemize}
\item The group {\it GT} of (structure-preserving) automorphisms of $B$
  acts on $A$ by composition on the right.
\item The group {\it GRT} of (structure-preserving) automorphisms of $C$    
  acts on $A$ by composition on the left.
\end{itemize}
\end{principle}

We apply this principle to a certain ``upgrade'' of the Kohno
isomorphism~\cite{Kohno:MonRep} (see also~\cite{KasselTuraev:Graphs})
between the unipotent completion $\widehat{PB}_n$ of the pure braid
group on $n$ strands and its associated graded algebra, which is a
certain completed algebra $\Apbc_n$ generated by symbols $t^{ij}$
with $1\leq i\neq j\leq n$.\footnote{In the language of Vassiliev
invariants, the Kohno isomorphism is a combination of three facts:
that the space of Vassiliev invariants of pure braids is the dual
of $\widehat{PB}_n$, that the associated graded space of Vassiliev
invariants of pure braids is dual to the algebra $\Apb_n$ of ``chord
diagrams'', and that the maps $\widehat{PB}_n\to\Apbc_n$ that we consider
are ``universal Vassiliev invariants''.} More specifically, in our case,
$B$ will be a certain category $\PB$ (defined in section~\ref{PBdef})
of {\em parenthesized braids}, and $C$ will be a certain category
$\PCD$ (defined in section~\ref{PCDdef}) of {\em parenthesized chord
diagrams}. On top of the category structure, both $\PB$ and $\PCD$ are
``fibered linear'', have natural ``basepoints'' (some specific morphisms
between some specific objects), natural ``coproducts'', and
natural ``extension'', ``cabling'', and ``strand removal''
operations, all defined in section~\ref{definitions}. Furthermore,
$\PB$ has a natural ``filtration'', $\PCD$ has a natural ``gradation''
(which induces a filtration as well), and these filtrations/gradations
(also defined in section~\ref{definitions}) respect all other structure
on $\PB$ and $\PCD$. In applying Principle~\ref{BasicPrinciple}, we will
only consider isomorphisms/automorphisms that respect all the additional
structure on $\PB$ and $\PCD$.

To be fair, we apply Principle~\ref{BasicPrinciple}
not to $B=\PB$ and $C=\PCD$, but rather
to their ``quotients'' $\PBm=\PB/\F{m+1}\PB$ and
$\PCDm=\PCD/\F{m+1}\PCD$ by their respective filtrations, or to
their ``completions'' $\PBc=\varprojlim_{m\to\infty}\PBm$ and
$\PCDc=\varprojlim_{m\to\infty}\PCDm$.  In section~\ref{ASS}
we show that every isomorphism (invertible structure-preserving
functor) $\Zc:\PBc\rightarrow\PCDc$ is determined by its action
on some specific morphism $a$ in $\PBc$, and that $\Zc(a)$ can be
interpreted as an associator.  We will thus identify the set of
all such $\Zc$'s with $\ASSc$, and get the two groups $\GTc$ and
$\GRTc$ (as well as their simple, transitive, and commuting actions)
entirely for free from Principle~\ref{BasicPrinciple}.  Similarly, using
Principle~\ref{BasicPrinciple} with $B=\PBm$ and $C=\PCDm$, we get groups
$\GTm$ and $\GRTm$ that act on the set $\ASSm$ of all ``associators up
to degree $m$''.

In section~\ref{Main} we start by explaining why the surjectivity of
the natural map $\pi:\GRTm\to\GRTmm$ implies the surjectivity of the map
$\ASSm\to\ASSmm$, which implies that there exists an iterative procedure
for finding an associator, and that a rational associator exists.

We then turn to the proof of the surjectivity of $\pi$. To do this, we
first write the relations defining $\GRTc$ explicitly. These turn out
to be the ``pentagon'', the ``classical hexagon'', the ``semi-classical
hexagon'', and some technical relations of lesser interest. It turns out
that the only relation that could challenge the surjectivity of $\pi$ is
the semi-classical hexagon, and so we spend the rest of section~\ref{Main}
proving that the semi-classical hexagon follows from the classical
hexagon, the pentagon, and the lesser relations. This is done by using
a certain 12-face polyhedron to show that the failure $\psi$ of the
semi-classical hexagon to hold lies in the kernel of some differential,
and by studying the relevant cohomology of the corresponding complex.

Just for completeness, in section~\ref{Just} we display the defining
formulas of $\GTc$ and $\GRTc$ that are not needed in the main argument. A
future part II of this paper will contain some additional results,
following~\cite[section 6]{Drinfeld:GalQQ}.

It is worthwhile to note that all our arguments depend on the existence of
at least one associator. Otherwise, we do not know that $\PBc$ and $\PCDc$
are at all isomorphic. So in a sense, all that we do is to take the
Knizhnik-Zamolodchikov associator $\Phi_{\KZ}$ (constructed by Drinfel'd)
and ``improve'' it.

Almost everything that we do appears either explicitly or implicitly in
Drinfel'd's paper~\cite{Drinfeld:GalQQ}. The presentation of $\GTc$  as a
group of automorphisms of some braid-group-like objects is due to Lochak
and Schneps~\cite{LochakSchneps:Automorphisms, LochakSchneps:Tower} (who
work with a different completion than ours).

\subsection{Acknowledgement} This paper grew out of a course I gave at
Harvard University in the spring semester of 1995, titled ``Knot Theory as
an Excuse''. One of the advertised goals of that course was to ``attempt
to read together two papers by Drinfel'd [\cite{Drinfeld:QuasiHopf,
Drinfeld:GalQQ}]''\footnote{All quotes taken from the official course
description.}, where I admitted that ``I have read about 20\% of the
material in these papers, understood about 20\% of what I read, and
got a lot out of it''. The idea was then to have ``a discussion group
in which everybody holds copies of the papers and we jointly try to
understand them''. Courses like that are usually doomed to fail, but
due to the amazing group of participants I think we managed to meet the
target of ``get ourself up to about 50\% on both figures [of reading
and understanding]''. These participants were:
D.~D.~Ben-Zvi,	
R.~Bott,	
A.~D'Andrea,	%
S.~Garoufalidis,
D.~J.~Goldberg, 
E.~Haley,	
M.~Hutchings,	
D.~Kazhdan,	
A.~Kirillov,	
T.~Kubo,	
S.~Majid,	
A.~Polishchuk,	
S.~Sternberg,	
D.~P.~Thurston, 
and
H.~L.~Wolfgang.	
I wish to thank them all for the part they took in the joint effort
that led to this paper.  In addition, I'd like to thank P.~Deligne,
E.~de-Shalit, and E.~Goren for teaching me some basic facts about
algebraic groups, and A.~Haviv, Yael~K., A.~Referee, E.~Rips, and
J.~D.~Stasheff for many useful comments.

%% file: definitions.tex
\section{The basic definitions} \lbl{definitions}

In this section we introduce the two mathematical structures $\PB$ and
$\PCD$ on which we will apply Principle~\ref{BasicPrinciple}. Let $A$ be
some fixed commutative associative ${\Bbb Q}$-algebra with unit (typically
${\Bbb C}$ or ${\Bbb Q}$). Most objects that we will define below ``have
coefficients'' in $A$. We will mostly suppress $A$ from the notation,
except in the few places where it matters.

\subsection{Parenthesized braids and $\protect\GT$} \lbl{PBdef}

A {\em parenthesized braid} is a braid (whose ends are points ordered
along a line) together with a parenthesization of its bottom end (the
{\em domain}) and its top end (the {\em range}). A {\em parenthesization}
of a sequence of points is a specification of a way of ``multiplying'' them
as if they were elements in a non-associative algebra. Rather then giving a
formal definition, Figure~\ref{ParenthesizedBraids} contains some
examples.

\begin{figure}[htpb]
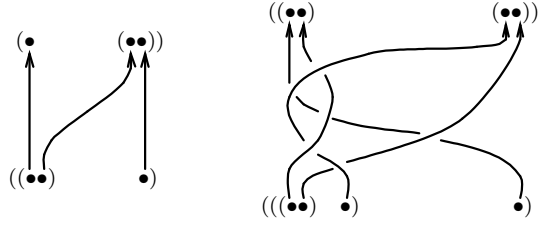

\def\b{\bullet}
\[ \printname{ParenthesizedBraids}
	\setlength{\unitlength}{0.6\standardunitlength}
	\begin{array}{c}  \hspace{-1.7mm}
        	\raisebox{-8pt}{\input draws/ParenthesizedBraids.tex }
        	\hspace{-1.9mm}
	\end{array}
 \]
\caption{A parenthesized braid whose domain is $((\b\b)\b)$ and whose
range is $(\b(\b\b))$ (left), and a parenthesized braid whose domain
is $(((\b\b)\b)\b)$ and whose range is $((\b\b)(\b\b))$ (right). Notice
that by convention we draw ``inner multiplications'' as closer endpoints,
and ``outer multiplications'' as farther endpoints. Below we will not
bother to specify the parenthesizations at the ends explicitly, as this
information can be read from the distance scales appearing in the way
we draw the ends.}
\lbl{ParenthesizedBraids}
\end{figure}

\par\noindent\par\noindent
\parbox{3.9in}{\setlength{\parindent}{\globalparindent}
  Parenthesized braids form a category in an obvious way. The objects of
  this category are parenthesizations, the morphisms are the parenthesized
  braids themselves, and composition is the operation of putting two
  parenthesized braid on top of each other, as on the right (provided the
  range of the first is the domain of the second).
} \qquad $
  {\def\arg{B_1}\printname{BoxMorphism}
	\setlength{\unitlength}{0.5\standardunitlength}
	\begin{array}{c}  \hspace{-1.7mm}
        	\raisebox{-8pt}{\input draws/BoxMorphism.tex }
        	\hspace{-1.9mm}
	\end{array}
}\!\!\circ
  {\def\arg{B_2}\printname{BoxMorphism}
	\setlength{\unitlength}{0.5\standardunitlength}
	\begin{array}{c}  \hspace{-1.7mm}
        	\raisebox{-8pt}{\input draws/BoxMorphism.tex }
        	\hspace{-1.9mm}
	\end{array}
} =
  \begin{array}{c}
    {\def\arg{B_2}\printname{BoxMorphism}
	\setlength{\unitlength}{0.5\standardunitlength}
	\begin{array}{c}  \hspace{-1.7mm}
        	\raisebox{-8pt}{\input draws/BoxMorphism.tex }
        	\hspace{-1.9mm}
	\end{array}
} \\
    \vspace{-9mm} \\
    {\def\arg{B_1}\printname{BoxMorphism}
	\setlength{\unitlength}{0.5\standardunitlength}
	\begin{array}{c}  \hspace{-1.7mm}
        	\raisebox{-8pt}{\input draws/BoxMorphism.tex }
        	\hspace{-1.9mm}
	\end{array}
}
  \end{array}
$

Furthermore, there are some naturally defined operations on
parenthesized braids. If $B$ is such a braid with $n$ strands, these
operations are:
\begin{itemize}
\item {\em Extension operations:} Let $d_0B=d^n_0B$ ($d_{n+1}B=d^n_{n+1}
B$) be $B$ with one straight strand added on the left (right), with ends
regarded as outer-most:
\[ d_0\left(\,\printname{extend1}
	\setlength{\unitlength}{0.5\standardunitlength}
	\begin{array}{c}  \hspace{-1.7mm}
        	\raisebox{-8pt}{\input draws/extend1.tex }
        	\hspace{-1.9mm}
	\end{array}
\right)=\printname{extend2}
	\setlength{\unitlength}{0.5\standardunitlength}
	\begin{array}{c}  \hspace{-1.7mm}
        	\raisebox{-8pt}{\input draws/extend2.tex }
        	\hspace{-1.9mm}
	\end{array}

  \quad;\qquad
  d_3\left(\,\printname{extend1}
	\setlength{\unitlength}{0.5\standardunitlength}
	\begin{array}{c}  \hspace{-1.7mm}
        	\raisebox{-8pt}{\input draws/extend1.tex }
        	\hspace{-1.9mm}
	\end{array}
\right)=\printname{extend3}
	\setlength{\unitlength}{0.5\standardunitlength}
	\begin{array}{c}  \hspace{-1.7mm}
        	\raisebox{-8pt}{\input draws/extend3.tex }
        	\hspace{-1.9mm}
	\end{array}
.
\]
\item {\em Cabling operations:} Let $d_iB=d^n_iB$ for $1\leq i\leq n$ be
the parenthesized braid obtained from $B$ by doubling its $i$th strand
(counting at the bottom), taking the ends of the resulting ``daughter
strands'' as an inner-most product:
\[ d_2\left(\,\printname{doubling1}
	\setlength{\unitlength}{0.5\standardunitlength}
	\begin{array}{c}  \hspace{-1.7mm}
        	\raisebox{-8pt}{\input draws/doubling1.tex }
        	\hspace{-1.9mm}
	\end{array}
\right)=\printname{doubling2}
	\setlength{\unitlength}{0.5\standardunitlength}
	\begin{array}{c}  \hspace{-1.7mm}
        	\raisebox{-8pt}{\input draws/doubling2.tex }
        	\hspace{-1.9mm}
	\end{array}
. \]
\item {\em Strand removal operations\footnote{The strand removal
operations (and all other $s_i$'s below) are important in the
applications, but play no crucial role in this paper and can be
systematically removed with no change to the end results.}:} Let
$s_iB=s^n_iB$ for $1\leq i\leq n$ be the parenthesized braid obtained
from $B$ by removing its $i$th strand (counting at the bottom):
\[ s_2\left(\,\printname{skeleton1}
	\setlength{\unitlength}{0.5\standardunitlength}
	\begin{array}{c}  \hspace{-1.7mm}
        	\raisebox{-8pt}{\input draws/skeleton1.tex }
        	\hspace{-1.9mm}
	\end{array}
\right)=\printname{removal2}
	\setlength{\unitlength}{0.5\standardunitlength}
	\begin{array}{c}  \hspace{-1.7mm}
        	\raisebox{-8pt}{\input draws/removal2.tex }
        	\hspace{-1.9mm}
	\end{array}
. \]
\end{itemize}

The {\em skeleton} ${\bold S}B$ of a parenthesized braid $B$ is the map
that it induces from the points of its domain to the points of its range,
taken together with the domain and range:
\begin{equation} \lbl{skeleton}
   {\bold S}\left(\,\printname{skeleton1}
	\setlength{\unitlength}{0.5\standardunitlength}
	\begin{array}{c}  \hspace{-1.7mm}
        	\raisebox{-8pt}{\input draws/skeleton1.tex }
        	\hspace{-1.9mm}
	\end{array}
\right) = \printname{skeleton2}
	\setlength{\unitlength}{0.5\standardunitlength}
	\begin{array}{c}  \hspace{-1.7mm}
        	\raisebox{-8pt}{\input draws/skeleton2.tex }
        	\hspace{-1.9mm}
	\end{array}
.
\end{equation}
More precisely, the skeleton ${\bold S}$ is a functor on the category of
parenthesized braids whose image is in the category $\PP$ of parenthesized
permutations, whose definition should be clear from its name and a simple
inspection of the example in \eqref{skeleton}. There are naturally defined
operations $d_i$ and $s_i$ on $\PP$ as in the case of parenthesized
braids, and the skeleton functor ${\bold S}$ intertwines the $d_i$'s
and the $s_i$'s acting on parenthesized braids and on parenthesized
permutations.

The category that we really need is a category of formal linear
combinations of parenthesized braids sharing the same skeleton:

\begin{definition} Let $\PB(A)=\PB$ (for {\bf Pa}renthesized {\bf B}raids)
be the category whose objects are parenthesizations and whose morphisms
are pairs $(P,\sum_{j=1}^k \beta_j B_j)$, where $P$ is a morphism in the
category of parenthesized permutations, the $B_j$'s are parenthesized
braids whose skeleton is $P$, and the $\beta_j$'s are coefficients in
the ground algebra $A$. The composition law in $\PB$ is the bilinear
extension of the composition law of parenthesized braids. There is a
natural forgetful ``skeleton'' functor ${\bold S}:\PB\to\PP$. If the sum
$\sum \beta_j B_j$ is not the empty sum, we usually suppress $P$ from the
notation, as it can be inferred from the $B_j$'s. See Figure~\ref{LinComb}.
\end{definition}

\begin{figure}[htpb]
\[ \printname{LinComb}
	\setlength{\unitlength}{0.75\standardunitlength}
	\begin{array}{c}  \hspace{-1.7mm}
        	\raisebox{-8pt}{\input draws/LinComb.tex }
        	\hspace{-1.9mm}
	\end{array}
 \]
\caption{A morphism $B$ in $\protect\PB$ and its skeleton ${\bold S}(B)$
  in $\protect\PP$.}
\lbl{LinComb}
\end{figure}

\subsubsection{Fibered linear categories}

The category $\PB$ together with the functor ${\bold S}:\PB\to\PP$ is
an example of a fibered linear category. Let ${\bold P}$ be a category
``of skeletons''. A {\em fibered linear category over ${\bold P}$} is a
pair $({\bold B},{\bold S}:{\bold B}\to{\bold P})$ of the form (category,
functor into ${\bold P}$), in which ${\bold B}$ has the same objects as
${\bold P}$, the ``skeleton'' functor ${\bold S}$ is the identity on
objects, the inverse image ${\bold S}^{-1}(P)$ of every morphism $P$ in
${\bold P}$ is a linear space, and so the composition maps in ${\bold
B}$ are bilinear in the natural sense. Many notions from the theory of
algebras have analogs for fibered linear categories, with the composition
of morphisms replacing the multiplication of elements. Let us list the
few such notions that we will use, without giving precise definitions:

\begin{myitemize}
\item A {\em subcategory} of a fibered linear category $({\bold B},{\bold
  S}:{\bold B}\to{\bold P})$ is a choice of a linear subspace in each
  ``space of morphisms with a fixed skeleton'' ${\bold S}^{-1}(P)$, so
  that the system of subspaces thus chosen is closed under composition.
\item An {\em ideal} in $({\bold B},{\bold S}:{\bold B}\to{\bold P})$ is a
  subcategory ${\bold I}$ so that if at least one of the two composable
  morphisms $B_1$ and $B_2$ in ${\bold B}$ is actually in ${\bold I}$, then
  the composition $B_1\circ B_2$ is also in ${\bold I}$.
\item One can take {\em powers} of ideals --- The morphisms of ${\bold
  I}^m$ will be all the morphisms in ${\bold B}$ that can be presented
  as compositions of $m$ morphisms in ${\bold I}$. The power ${\bold
  I}^m$ is also an ideal in ${\bold B}$.
\item One can form the {\em quotient} ${\bold B}/{\bold I}$ of a fibered
  linear category ${\bold B}$ by an ideal ${\bold I}$ in it, and the
  result is again a fibered linear category.
\item {\em Direct sums} of fibered linear categories that are fibered over
  the same skeleton category can be formed.
\item One can define {\em filtered} and {\em graded} fibered linear
  categories. One can talk about the {\em associated graded} fibered linear
  category of a given filtered fibered linear category.
\item One can take the {\em inverse limit} of an inverse system of
  fibered linear categories (fibered in a compatible way over the same
  category of skeletons). In particular, if ${\bold I}$ is an ideal
  in a fibered linear category ${\bold B}$, one can form ``the ${\bold
  I}$-adic completion $\hat{{\bold B}}=\varprojlim_{m\to\infty}{\bold
  B}/{\bold I}^m$. The ${\bold I}$-adic completion is a filtered fibered
  linear category.
\item {\em Tensor powers} of a fibered linear category $({\bold
  B},{\bold S}:{\bold B}\to{\bold P})$ can be defined. For example,
  ${\bold B}\otimes{\bold B}$ will have the same set of objects as
  ${\bold B}$, and for any two such objects $O_1$ and $O_2$, we set
  \[ \mor_{{\bold B}\otimes{\bold B}}(O_1,O_2)=
    \coprod_{P\in\mor_{{\bold P}}(O_1,O_2)}
    {\bold S}^{-1}(P)\otimes{\bold S}^{-1}(P).
  \]
  ${\bold B}\otimes{\bold B}$ is again a fibered linear category.
\item The notion of a {\em coproduct functor}
  ${\bold \Box}:{\bold B}\to{\bold B}\otimes{\bold B}$ makes sense.
\end{myitemize}

\subsubsection{Back to parenthesized braids}

We can now introduce some more structure on $\PB$, and specify completely
the mathematical structures that will play the role of $B$ in
Principle~\ref{BasicPrinciple}.

\begin{definition} Let ${\bold \Box}:\PB\to\PB\otimes\PB$ be the coproduct
functor defined by setting each individual parenthesized braid $B$ to
be {\em group-like}, that is, by setting ${\bold \Box}(B)=B\otimes B$.
\end{definition}

Let ${\bold I}$ be the {\em augmentation ideal} of $\PB$, the ideal of
all pairs $(P,\sum\beta_j B_j)$ in which $\sum\beta_j=0$. Powers of this
ideal define the {\em unipotent filtration} of $\PB$: $\F{m}\PB={\bold
I}^{m+1}$.

\begin{definition} Let $\PBm=\PB/\F{m}\PB=\PB/{\bold I}^{m+1}$ be the $m$th
{\em unipotent quotient}\footnote{If you are familiar with Vassiliev
invariants, notice that $\PBm$ is simply $\PB$ moded out by
``$(m+1)$-singular parenthesized braids''.} of $\PB$, and let
$\PBc=\varprojlim_{m\to\infty}\PBm$ be the {\em unipotent completion}
of $\PB$.
\end{definition}

Let $\sigma$ be the parenthesized braid $\printname{sigma}
	\setlength{\unitlength}{0.25\standardunitlength}
	\begin{array}{c}  \hspace{-1.7mm}
        	\raisebox{-8pt}{\input draws/sigma.tex }
        	\hspace{-1.9mm}
	\end{array}
$.

The fibered linear categories $\PBm$ and $\PBc$ inherit the operations
$d_i$ and $s_i$ from parenthesized braids, and a coproduct ${\bold
\Box}$ and a filtration $\F{\star}$ from $\PB$. The specific parenthesized
braid $\sigma$ can be regarded as a morphism in any of these categories.

\begin{definition} Let $\GTm$ and $\GTc$ (really, $\GTm(A)$ and $\GTc(A)$)
be the groups of structure preserving automorphisms of $\PBm$ and
$\PBc$, respectively. That is, the groups of all functors $\PBm\to\PBm$
(or $\PBc\to\PBc$) that cover the skeleton functor, intertwine $d_i$,
$s_i$ and ${\bold \Box}$ and fix $\sigma$. In short, let
\begin{eqnarray*}
  B^{(m)} & = & \left(
    \PBm,{\bold S}:\PBm\to\PP,d_i,s_i,{\bold \Box},\sigma \right); \\
  \hat{B} & = & \left(
    \PBc,{\bold S}:\PBc\to\PP,d_i,s_i,{\bold \Box},\sigma \right);
\end{eqnarray*}
\[ \GTm=\Aut B^{(m)}; \qquad \GTc=\Aut\hat{B}. \]
\end{definition}

\begin{remark} One easily sees that elements of $\GTm$ ($\GTc$)
automatically preserve the filtration $\F{\star}$.
\end{remark}

\begin{claim} $\PB$ is generated by $a^{\pm 1}$, $\sigma^{\pm 1}$, and
their various images by repeated applications of the $d_i$'s, where
\[ a=\printname{a}
	\setlength{\unitlength}{0.5\standardunitlength}
	\begin{array}{c}  \hspace{-1.7mm}
        	\raisebox{-8pt}{\input draws/a.tex }
        	\hspace{-1.9mm}
	\end{array}
, \qquad \sigma=\printname{sigma}
	\setlength{\unitlength}{0.5\standardunitlength}
	\begin{array}{c}  \hspace{-1.7mm}
        	\raisebox{-8pt}{\input draws/sigma.tex }
        	\hspace{-1.9mm}
	\end{array}
. \]
\end{claim}

\begin{pf} (sketch) The main point is that any of the standard generators
of the braid group can be written in terms of $a^{\pm 1}$ and $\sigma^{\pm
1}$ and their images. For example,
\[ \printname{BraidGenerator}
	\setlength{\unitlength}{0.5\standardunitlength}
	\begin{array}{c}  \hspace{-1.7mm}
        	\raisebox{-8pt}{\input draws/BraidGenerator.tex }
        	\hspace{-1.9mm}
	\end{array}
=\printname{BGRewritten}
	\setlength{\unitlength}{0.5\standardunitlength}
	\begin{array}{c}  \hspace{-1.7mm}
        	\raisebox{-8pt}{\input draws/BGRewritten.tex }
        	\hspace{-1.9mm}
	\end{array}

  = d_0a^{-1}\circ d_0d_3\sigma\circ d_0a.
\]
\vskip -1cm
\end{pf}
\vskip 1cm

\subsection{Parenthesized chord diagrams and $\protect\GRT$} \lbl{PCDdef}

The category $\PCD$, the main ingredient of the mathematical object
$C$ on which we will apply Principle~\ref{BasicPrinciple}, can be
viewed as natural in two (equivalent) ways. First, $\PCD$ is natural
because it is the associated graded of $\PB$, as will be proven in
section~\ref{ASS}.  $\PCD$ can also be viewed as the category of ``chord
diagrams for finite-type (Vassiliev) invariants \cite{Bar-Natan:Vassiliev,
Bar-Natan:Braids, Birman:Bulletin, BirmanLin:Vassiliev, Goussarov:New,
Goussarov:nEquivalence, Kontsevich:Vassiliev, Vassiliev:CohKnot,
Vassiliev:Book} of parenthesized braids'', and all the operations that we
will define on $\PCD$ are inherited from their parallels on parenthesized
braids, that were defined in section~\ref{PBdef}. I prefer not to make
more than a few comments about the latter viewpoint below. Saying more
requires repeating well known facts about finite-type invariants, and
these can easily be found in the literature.  If you already know about
Vassiliev invariants and chord diagrams, you'll find the relation between
them and the definitions below rather clear.  Unfortunately, if
finite-type invariants are not mentioned, we have to start with some
unmotivated definitions.

\begin{definition} \lbl{def:Apb}
Let $\Apb_n=\Apb_n(A)$ be the algebra (over the ground algebra $A$)
generated by symbols $t^{ij}$ for $1\leq i\neq j\leq n$, subject to the
relations $t^{ij}=t^{ji}$, $[t^{ij},t^{kl}]=0$ if $|\{i,j,k,l\}|=4$,
and $[t^{jk},t^{ij}+t^{ik}]=0$ if $|\{i,j,k\}|=3$. The algebra
$\Apb_n$ is graded by setting $\deg t^{ij}=1$; let $\G{m}\Apb_n$
be the degree $m$ piece of $\Apb_n$, let $\F{\star}\Apb_n$ be the
filtration defined by $\F{m}\Apb_n=\bigoplus_{m'>m}\G{m'}\Apb_n$,
let $\Apbm_n$ be $\Apb_n/\F{m}\Apb_n$, and let $\Apbc_n$ be the
graded completion $\varprojlim_{m\to\infty}\Apbm_n$ of $\Apb_n$. We
call elements of $\Apb_n$ {\em chord diagrams}, and draw them as in
Figure~\ref{ChordDiagram}. (In the language of finite-type invariants,
$\Apb_n$ is the algebra of chord diagrams for $n$-strand pure braids,
and the last relation is the ``$4T$'' relation.)
\end{definition}

\begin{figure}[htpb]
\[ t^{13}t^{13}t^{12}t^{23} \quad\longleftrightarrow\quad
  \printname{ChordDiagram}
	\setlength{\unitlength}{0.5\standardunitlength}
	\begin{array}{c}  \hspace{-1.7mm}
        	\raisebox{-8pt}{\input draws/ChordDiagram.tex }
        	\hspace{-1.9mm}
	\end{array}
; \qquad
  4T:\quad\printname{4T}
	\setlength{\unitlength}{0.5\standardunitlength}
	\begin{array}{c}  \hspace{-1.7mm}
        	\raisebox{-8pt}{\input draws/4T.tex }
        	\hspace{-1.9mm}
	\end{array}

\]
\caption{Elements of $\Apb_3$ are presented as chord diagrams made of $3$
vertical strands and some number of horizontal chords connecting them. A
chord connecting the $i$th strand to the $j$th strand represents $t^{ij}$,
and products are read from the bottom to the top of the diagram.}
\lbl{ChordDiagram}
\end{figure}

\begin{definition} \lbl{PermutationAction} There is an action of the
symmetric group ${\cal S}_n$ on $\Apb_n$ by ``permuting the vertical
strands'', denoted by $(\tau,\Psi)\mapsto\Psi^\tau$:
\[ \Psi=\printname{Psi}
	\setlength{\unitlength}{0.5\standardunitlength}
	\begin{array}{c}  \hspace{-1.7mm}
        	\raisebox{-8pt}{\input draws/Psi.tex }
        	\hspace{-1.9mm}
	\end{array}
\mapsto\Psi^{231}=\printname{Psi231}
	\setlength{\unitlength}{0.5\standardunitlength}
	\begin{array}{c}  \hspace{-1.7mm}
        	\raisebox{-8pt}{\input draws/Psi231.tex }
        	\hspace{-1.9mm}
	\end{array}
. \]
\end{definition}

\begin{definition} \lbl{Action} Let $d_i=d^n_i:\Apb_n\to\Apb_{n+1}$
for $0\leq i\leq n+1$ and $s_i=s^n_i:\Apb_n\to\Apb_{n-1}$ for $1\leq i\leq
n$ be the algebra morphisms defined by their action on the generators
$t^{jk}$ (with $j<k$) as follows:
\[ d_it^{jk}=\begin{cases}
    t^{j+1,k+1}			& i<j<k \\
    t^{j,k+1}+t^{j+1,k+1}	& i=j<k \\
    t^{j,k+1}			& j<i<k \\
    t^{jk}+t^{j,k+1}		& j<i=k \\
    t^{jk}			& j<k<i
  \end{cases}\qquad
  s_it^{jk}=\begin{cases}
    t^{j-1,k-1}	& i<j<k \\
    0		& i=j<k \\
    t^{j,k-1}	& j<i<k \\
    0		& j<i=k \\
    t^{jk}	& j<k<i.
  \end{cases}
\]
Graphically, $d^n_0$ ($d^n_{n+1}$) acts by adding a strand on the left
(right), $d^n_i$ for $1\leq i\leq n$ acts by doubling the $i$th strand and
summing all the possible ways of lifting the chords that were connected to
the $i$th strand to the two daughter strands, and $s^n_i$ acts by deleting
the $i$th strand and mapping the chord diagram to $0$ if any chord in it
was connected to the $i$th strand:
\[ d_0\left(\,\printname{Action1}
	\setlength{\unitlength}{0.5\standardunitlength}
	\begin{array}{c}  \hspace{-1.7mm}
        	\raisebox{-8pt}{\input draws/Action1.tex }
        	\hspace{-1.9mm}
	\end{array}
\right)=\printname{Action2}
	\setlength{\unitlength}{0.5\standardunitlength}
	\begin{array}{c}  \hspace{-1.7mm}
        	\raisebox{-8pt}{\input draws/Action2.tex }
        	\hspace{-1.9mm}
	\end{array}
;
  \qquad d_2\left(\,\printname{Action1}
	\setlength{\unitlength}{0.5\standardunitlength}
	\begin{array}{c}  \hspace{-1.7mm}
        	\raisebox{-8pt}{\input draws/Action1.tex }
        	\hspace{-1.9mm}
	\end{array}
\right)=
  \printname{Action4}
	\setlength{\unitlength}{0.5\standardunitlength}
	\begin{array}{c}  \hspace{-1.7mm}
        	\raisebox{-8pt}{\input draws/Action4.tex }
        	\hspace{-1.9mm}
	\end{array}
=\printname{Action3}
	\setlength{\unitlength}{0.5\standardunitlength}
	\begin{array}{c}  \hspace{-1.7mm}
        	\raisebox{-8pt}{\input draws/Action3.tex }
        	\hspace{-1.9mm}
	\end{array}
;
\]
\[ s_1\left(\,\printname{Action2}
	\setlength{\unitlength}{0.5\standardunitlength}
	\begin{array}{c}  \hspace{-1.7mm}
        	\raisebox{-8pt}{\input draws/Action2.tex }
        	\hspace{-1.9mm}
	\end{array}
\right)=\printname{Action1}
	\setlength{\unitlength}{0.5\standardunitlength}
	\begin{array}{c}  \hspace{-1.7mm}
        	\raisebox{-8pt}{\input draws/Action1.tex }
        	\hspace{-1.9mm}
	\end{array}
; \qquad
  s_1\left(\,\printname{Action1}
	\setlength{\unitlength}{0.5\standardunitlength}
	\begin{array}{c}  \hspace{-1.7mm}
        	\raisebox{-8pt}{\input draws/Action1.tex }
        	\hspace{-1.9mm}
	\end{array}
\right)=0.
\]
(Here and below the symbol $\setlength{\unitlength}{0.5\standardunitlength}
	\begin{array}{c}  \hspace{-1.7mm}
        	\raisebox{-2pt}{\input draws/BoxSymbol.tex }
        	\hspace{-1.9mm}
	\end{array}
$ means
$\setlength{\unitlength}{0.5\standardunitlength}
	\begin{array}{c}  \hspace{-1.7mm}
        	\raisebox{-2pt}{\input draws/BoxMeans.tex }
        	\hspace{-1.9mm}
	\end{array}
$).
\end{definition}

\begin{definition} Let $\Box:\Apb_n\to\Apb_n\otimes\Apb_n$ be the
coproduct defined by declaring the $t^{ij}$'s to be primitive:
$\Box(t^{ij})=t^{ij}\otimes 1 + 1\otimes t^{ij}$.
\end{definition}

\begin{definition} $\PCD=\PCD(A)$ (for {\bf Pa}renthesized
{\bf C}hord {\bf D}iagrams) is the category whose objects are
parenthesizations and whose morphisms are formal products $D\cdot
P$, where $P$ is a parenthesized permutation of $n$ objects (for
some $n$) and $D\in\Apb_n(A)$. The composition law in $\PCD$ is
$D_1\cdot P_1\circ D_2\cdot P_2=(D_1\cdot D_2^{P_1})\cdot(P_1\circ
P_2)$ (whenever $P_1$ and $P_2$ are composable), where $D_2^{P_1}$
denotes the action of of the permutation $P_1$ on $D_2$ as in
Definition~\ref{PermutationAction}. This composition law is better
seen graphically as in Figure~\ref{CompositionLaw}. $\PCD$ inherits
a grading $\PCD=\bigoplus_{m}\G{m}\PCD$ from $\Apb_\star$, and
is fibered linear over $\PP$ with the skeleton functor ${\bold
S}:D\cdot P\mapsto P$. $\PCD$ is also be filtered by setting
$\F{m}\PCD=\bigoplus_{m'>m}\G{m'}\PCD$. $\PCD$ inherits a coproduct
${\bold \Box}:\PCD\to\PCD\otimes\PCD$ from the coproduct $\Box$ of
$\Apb_n$.
\end{definition}

\begin{figure}[htpb]
\[
  \left(t^{12}\cdot\setlength{\unitlength}{0.5\standardunitlength}
	\begin{array}{c}  \hspace{-1.7mm}
        	\raisebox{-2pt}{\input draws/Comp1.tex }
        	\hspace{-1.9mm}
	\end{array}
\right)\circ
    \left(t^{23}t^{13}\cdot\setlength{\unitlength}{0.5\standardunitlength}
	\begin{array}{c}  \hspace{-1.7mm}
        	\raisebox{-2pt}{\input draws/Comp2.tex }
        	\hspace{-1.9mm}
	\end{array}
\right)
  \to \setlength{\unitlength}{0.5\standardunitlength}
	\begin{array}{c}  \hspace{-1.7mm}
        	\raisebox{-2pt}{\input draws/Comp3.tex }
        	\hspace{-1.9mm}
	\end{array}
\circ\setlength{\unitlength}{0.5\standardunitlength}
	\begin{array}{c}  \hspace{-1.7mm}
        	\raisebox{-2pt}{\input draws/Comp4.tex }
        	\hspace{-1.9mm}
	\end{array}

  \to \setlength{\unitlength}{0.5\standardunitlength}
	\begin{array}{c}  \hspace{-1.7mm}
        	\raisebox{-2pt}{\input draws/Comp5.tex }
        	\hspace{-1.9mm}
	\end{array}

  \to \setlength{\unitlength}{0.5\standardunitlength}
	\begin{array}{c}  \hspace{-1.7mm}
        	\raisebox{-2pt}{\input draws/Comp6.tex }
        	\hspace{-1.9mm}
	\end{array}

  \to t^{12}t^{23}t^{13}\cdot\setlength{\unitlength}{0.5\standardunitlength}
	\begin{array}{c}  \hspace{-1.7mm}
        	\raisebox{-2pt}{\input draws/Comp7.tex }
        	\hspace{-1.9mm}
	\end{array}

\]
\caption{The composition of a morphism in
  $\protect\mor_{\protect\PCD}
    ((\bullet(\bullet\bullet)),(\bullet(\bullet\bullet)))$
  with a morphism in
  $\protect\mor_{\protect\PCD}
    ((\bullet(\bullet\bullet)),((\bullet\bullet)\bullet))$.
} \lbl{CompositionLaw}
\end{figure}

\begin{definition} \lbl{PCDOperators}
As in the case of $\PB$, there are some naturally defined
operations on $\PCD$. If $D\cdot P$ is a parenthesized chord diagram on
$n$ strands, set $d_i(D\cdot P)=d^n_i(D\cdot P)=d^n_iD\cdot d^n_iP$,
and similarly for $s_i=s^n_i$. These operations are:
\begin{itemize}
\item {\em Extension operations:} $d_0$ ($d_{n+1}$) adds a far-away
  independent strand on the left (right).
\item {\em Cabling operations:} $d_i B$ with $1\leq i\leq n$ doubles
  the $i$th strand and sums all possible ways of lifting the chords that
  were connected to the $i$th strand to the two daughter strands.
\item {\em Strand removal operations:} $s_i$ removes the $i$th strand and
  maps everything to $0$ if there was any chord connected to the $i$th
  strand.
\end{itemize}
\end{definition}

\begin{definition} Let $\PCDm$ be the category $\PCD/\F{m}\PCD$ of 
parenthesized chord diagrams of degree up to $m$, and let $\PCDc$ be the
category $\varprojlim_{m\to\infty}\PCDm$ of formal power series of
parenthesized chord diagrams. The fibered linear categories $\PCDm$
and $\PCDc$ inherit the operations $d_i$ and $s_i$, the coproduct
${\bold \Box}$ and the filtration $\F{\star}$ from $\PCD$.
\end{definition}

Let $X$ and $H$ be the parenthesized chord diagrams $\printname{X}
	\setlength{\unitlength}{0.25\standardunitlength}
	\begin{array}{c}  \hspace{-1.7mm}
        	\raisebox{-8pt}{\input draws/X.tex }
        	\hspace{-1.9mm}
	\end{array}
$
and $\printname{H}
	\setlength{\unitlength}{0.25\standardunitlength}
	\begin{array}{c}  \hspace{-1.7mm}
        	\raisebox{-8pt}{\input draws/H.tex }
        	\hspace{-1.9mm}
	\end{array}
$ respectively, and let $\tilde{R}$ be the formal
exponential $\tilde{R}=\exp\left(\frac{1}{2}H\right)\cdot X$, regarded
a morphism in $\PCDm$ or $\PCDc$.

\begin{definition} Let $\GRTm$ and $\GRTc$ (really, $\GRTm(A)$ and
$\GRTc(A)$) be the groups of structure preserving automorphisms of
$\PCDm$ and $\PCDc$, respectively. That is, the groups of all functors
$\PCDm\to\PCDm$ (or $\PCDc\to\PCDc$) that cover the skeleton functor,
intertwine $d_i$, $s_i$ and ${\bold \Box}$ and fix $\tilde{R}$. In
short, let
\begin{eqnarray*} C^{(m)} & = & \left(
    \PCDm,{\bold S}:\PCDm\to\PP,d_i,s_i,{\bold \Box},\tilde{R}
  \right); \\
  \hat{C} & = & \left(
    \PCDc,{\bold S}:\PCDc\to\PP,d_i,s_i,{\bold \Box},\tilde{R}
  \right);
\end{eqnarray*}
\[ \GRTm=\Aut C^{(m)}; \qquad \GRTc=\Aut\hat{C}. \]
\end{definition}

\begin{remark} \lbl{FixXH} Elements of $\GRTm$ ($\GRTc$) fix each of $X$
and $H$ individually. Indeed, $\tilde{R}^2=\exp H$ and hence $\exp H$
and thus $H$ are fixed. But then $X=\exp(-\frac{1}{2}H)\tilde{R}$ is
fixed too.
\end{remark}

\begin{claim} \lbl{PaCDGens} $\PCD$ is generated by $a^{\pm 1}$, $X$, $H$,
and their various images by repeated applications of the $d_i$'s, where
\[ a=\printname{a}
	\setlength{\unitlength}{0.5\standardunitlength}
	\begin{array}{c}  \hspace{-1.7mm}
        	\raisebox{-8pt}{\input draws/a.tex }
        	\hspace{-1.9mm}
	\end{array}
, \qquad X=\printname{X}
	\setlength{\unitlength}{0.5\standardunitlength}
	\begin{array}{c}  \hspace{-1.7mm}
        	\raisebox{-8pt}{\input draws/X.tex }
        	\hspace{-1.9mm}
	\end{array}
, \qquad H=\printname{H}
	\setlength{\unitlength}{0.5\standardunitlength}
	\begin{array}{c}  \hspace{-1.7mm}
        	\raisebox{-8pt}{\input draws/H.tex }
        	\hspace{-1.9mm}
	\end{array}
. \]
(Notice that the symbol ``$a$'' plays a double role, as a generator of $\PB$
and as a generator of $\PCD$).
\end{claim}

\begin{pf} (sketch) Perhaps one illustrative example will suffice:
\[ \printname{t13}
	\setlength{\unitlength}{0.5\standardunitlength}
	\begin{array}{c}  \hspace{-1.7mm}
        	\raisebox{-8pt}{\input draws/t13.tex }
        	\hspace{-1.9mm}
	\end{array}
=\printname{t13Rewritten}
	\setlength{\unitlength}{0.5\standardunitlength}
	\begin{array}{c}  \hspace{-1.7mm}
        	\raisebox{-8pt}{\input draws/t13Rewritten.tex }
        	\hspace{-1.9mm}
	\end{array}

  = d_0X\circ a^{-1}\circ d_3H\circ a\circ d_0X.
\]
\vskip -8mm
\end{pf}

\begin{remark} Remark~\ref{FixXH} and claim~\ref{PaCDGens} imply that
elements of $\GRTm$ ($\GRTc$) automatically preserve the filtration
$\F{\star}$.
\end{remark}

%% file: draws/sigma.tex
%
\begingroup\makeatletter\ifx\SetFigFont\undefined%
\gdef\SetFigFont#1#2#3#4#5{%
  \reset@font\fontsize{#1}{#2pt}%
  \fontfamily{#3}\fontseries{#4}\fontshape{#5}%
  \selectfont}%
\fi\endgroup%
\begin{picture}(1224,1239)(0,-10)
\thicklines
\path(1212,12)(762,462)
\path(462,762)(12,1212)
\path(133.976,1143.057)(12.000,1212.000)(80.943,1090.024)
\path(1143.057,1090.024)(1212.000,1212.000)(1090.024,1143.057)
\path(1212,1212)(12,12)
\end{picture}

%% file: draws/X.tex
%
\begingroup\makeatletter\ifx\SetFigFont\undefined%
\gdef\SetFigFont#1#2#3#4#5{%
  \reset@font\fontsize{#1}{#2pt}%
  \fontfamily{#3}\fontseries{#4}\fontshape{#5}%
  \selectfont}%
\fi\endgroup%
\begin{picture}(1224,1239)(0,-10)
\thicklines
\path(12,12)(1212,1212)
\path(1143.057,1090.024)(1212.000,1212.000)(1090.024,1143.057)
\path(1212,12)(12,1212)
\path(133.976,1143.057)(12.000,1212.000)(80.943,1090.024)
\end{picture}

%% file: draws/H.tex
%
\begingroup\makeatletter\ifx\SetFigFont\undefined%
\gdef\SetFigFont#1#2#3#4#5{%
  \reset@font\fontsize{#1}{#2pt}%
  \fontfamily{#3}\fontseries{#4}\fontshape{#5}%
  \selectfont}%
\fi\endgroup%
\begin{picture}(1224,1239)(0,-10)
\thicklines
\path(12,612)(1212,612)
\path(12,12)(12,1212)
\path(49.500,1077.000)(12.000,1212.000)(-25.500,1077.000)
\path(1212,12)(1212,1212)
\path(1249.500,1077.000)(1212.000,1212.000)(1174.500,1077.000)
\end{picture}

%% file: ASS.tex
\section{Isomorphisms and associators} \lbl{ASS}

In this section we make the key observation that makes
Principle~\ref{BasicPrinciple} useful in our case: The fact that
the set of all associators \`a la Drinfel'd~\cite{Drinfeld:QuasiHopf,
Drinfeld:GalQQ} can be identified with the set of all structure-preserving
functors $\Zc:\hat{B}\to\hat{C}$. Recall that $A$ is some fixed commutative
associative ${\Bbb Q}$-algebra with unit.

\begin{definition} An {\em associator} is an invertible element $\Phi$
of $\Apbc_3(A)$ satisfying the following axioms:
\begin{itemize}
\item The {\em pentagon} axiom holds in $\Apbc_4$:
  \begin{equation}
    d_4\Phi\cdot d_2\Phi\cdot d_0\Phi
    = d_1\Phi\cdot d_3\Phi.
  \tag{$\pentagon$}
  \end{equation}
\item The {\em hexagon} axioms hold in $\Apbc_3$:
  \begin{equation}
    d_1\exp\left(\pm\frac{1}{2}t^{12}\right) =
    \Phi\cdot\exp\left(\pm\frac{1}{2}t^{23}\right)\cdot(\Phi^{-1})^{132}
    \cdot\exp\left(\pm\frac{1}{2}t^{13}\right)\cdot\Phi^{312}.
  \tag{$\hexagon\!_\pm$}
  \end{equation}
\item $\Phi$ is {\em non-degenerate}: $s_1\Phi=s_2\Phi=s_3\Phi=1$.
\item $\Phi$ is {\em group-like}: $\Box\Phi=\Phi\otimes\Phi$.
\end{itemize}
(Apart from the different notation, this definition is equivalent to
Drinfel'd's~\cite{Drinfeld:GalQQ} definition of an $Fr(A,B)$-valued
$\varphi$, and practically equivalent to the definition of an $\text{\bf
AP}^{\text{\it hor}}$-valued $\Phi$ in~\cite{Bar-Natan:NAT}.)
\end{definition}

\begin{definition} Let $\ASSc=\ASSc(A)$ be the set of associators
$\Phi\in\Apbc_3(A)$. Similarly, if we mod out by degrees higher than $m$,
we can define {\em associators up to degree $m$} and the set $\ASSm$.
\end{definition}

\begin{remark} The hexagon axiom for $\Phi\in\ASSc$ or $\Phi\in\ASSm$
implies that $\Phi=1+$(higher degree terms).
\end{remark}

By the definition of $\hat{B}$ and $\hat{C}$, a structure-preserving
functor $\Zc:\hat{B}\to\hat{C}$ carries $\sigma$ to $\tilde{R}$, and thus
it is determined by its value $\Zc(a)$ on the remaining generator of $\PB$.
As $\Zc$ must cover the skeleton functor, $\Zc(a)$ must be of the form
$\Phi_\Zc\cdot a$, for some $\Phi_\Zc\in\Apbc_3$.

\begin{proposition} \lbl{KeyObservation}
If $\Zc$ is a structure preserving functor $\hat{B}\to\hat{C}$, then
$\Phi_\Zc$ is an associator, and the map $\Zc\mapsto\Phi_\Zc$ is
a bijection between the set of all structure-preserving functors
$\Zc:\hat{B}\to\hat{C}$ and the set $\ASSc$ of all associators
$\Phi\in\Apbc_3$. A similar construction can be made in the case
of $B^{(m)}$, $C^{(m)}$ and $\ASSm$, and the same statements hold.
\end{proposition}

Before we can prove Proposition~\ref{KeyObservation}, we need a bit more
insight about the structure of $\Apb_n$.

\begin{lemma} \lbl{Locality} The following two relations hold in $\Apb_n$:
\begin{enumerate}
\item {\em Locality in space:} For any $k\leq n$, the subalgebra of
  $\Apb_n$ generated by $\{t^{ij}:i,j\leq k\}$ commutes with the subalgebra
  generated by $\{t^{ij}:i,j>k\}$. In pictures, we see that elements that
  live in ``different parts of space'' commute:
  \[ \printname{ApbLocSpace}
	\setlength{\unitlength}{0.5\standardunitlength}
	\begin{array}{c}  \hspace{-1.7mm}
        	\raisebox{-8pt}{\input draws/ApbLocSpace.tex }
        	\hspace{-1.9mm}
	\end{array}
 \]
\item {\em Locality in scale} Elements that live in ``different scales''
  commute. This is best explained by a picture, with notation as in
  Definition~\ref{Action}:
  \[ \printname{ApbLocScale}
	\setlength{\unitlength}{0.7\standardunitlength}
	\begin{array}{c}  \hspace{-1.7mm}
        	\raisebox{-8pt}{\input draws/ApbLocScale.tex }
        	\hspace{-1.9mm}
	\end{array}
 \]
  (We think of the part $A$ as ``local'', as it involves only the ``local''
  group of strands, and of the rest as ``global'', as it regards the
  ``local'' group of strands as ``equal''.)
\end{enumerate}
(A similar statement is \cite[Lemma 3.4]{Bar-Natan:NAT}.)
\end{lemma}

\begin{pf*}{Proof of Lemma~\ref{Locality}}
Locality in space follows from repeated application of the relation
$t^{ij}t^{kl}=t^{kl}t^{ij}$ with $i<j<k<l$. Locality in scale follows from
repeated application of the relation $t^{ij}t^{kl}=t^{kl}t^{ij}$ with
general $i,j,k,l$ with $|\{i,j,k,l\}|=4$, and the $4T$ relation, which
can be redrawn in the more suggestive form $\setlength{\unitlength}{0.6\standardunitlength}
	\begin{array}{c}  \hspace{-1.7mm}
        	\raisebox{-2pt}{\input draws/Suggestive4T.tex }
        	\hspace{-1.9mm}
	\end{array}
$.
\end{pf*}

\begin{pf*}{Proof of Proposition~\ref{KeyObservation}}
Let $\Zc$ be a structure preserving functor $\hat{B}\to\hat{C}$, and let
$\Phi=\Phi_\Zc$. Apply $\Zc$ to the parenthesized braid equality
\[ \printname{Pentagon}
	\setlength{\unitlength}{0.5\standardunitlength}
	\begin{array}{c}  \hspace{-1.7mm}
        	\raisebox{-8pt}{\input draws/Pentagon.tex }
        	\hspace{-1.9mm}
	\end{array}
\quad;\qquad
  d_4a\circ d_2a\circ d_0a = d_1a\circ d_3a,
\]
and, using $\Zc(a)=\Phi\cdot a$, get
\[ (d_4\Phi\cdot d_2\Phi\cdot d_0\Phi)
    \cdot(d_4a\circ d_2a\circ d_0a)
  =(d_1\Phi\cdot d_3\Phi)\cdot(d_1a\circ d_3a).
\]
The $\Apbc$ part of this equality is precisely the fact the $\Phi$
satisfies the pentagon equation.

Similarly, the parenthesized braid equality
\begin{equation} \lbl{HexagonPlus}
  \printname{HexagonPlus}
	\setlength{\unitlength}{0.5\standardunitlength}
	\begin{array}{c}  \hspace{-1.7mm}
        	\raisebox{-8pt}{\input draws/HexagonPlus.tex }
        	\hspace{-1.9mm}
	\end{array}
\quad;\qquad
  d_1\sigma =
    a\circ d_0\sigma\circ a^{-1}\circ d_3\sigma\circ a
\end{equation}
together with $\Zc(a)=\Phi\cdot a$ and $\Zc(\sigma)=\tilde{R}$ implies 
the positive hexagon equation $\hexagon\!_+$. Likewise, the same
parenthesized braid equality but with $\sigma$ replaced by $\sigma^{-1}$
implies $\hexagon\!_-$.

$s_1a=s_2a=s_3a$ is the identity morphism in
$\mor((\bullet\bullet),(\bullet\bullet))$, and after applying $\Zc$ we
find that $\Phi$ is non-degenerate. Finally, $a$ is group-like in $\PB$
and as $\Zc$ preserves the coproduct, $\Phi$ is also group-like. Hence
we have verified that $\Phi_\Zc=\Phi$ is an associator.

To show that the map $\Zc\mapsto\Phi_\Zc$ is a bijection we construct
an inverse map. Let $\Phi$ be an associator. We try to define a
functor $\Zc=\Zc_\Phi:\PBc\to\PCDc$ by setting $\Zc(a)=a\cdot\Phi$ and
$\Zc(\sigma)=\tilde{R}$, and by extending it to all other generators of
$\PB$ in a way compatible with the $d_i$'s. We need to verify
that this extension yields a well-defined functor; that is, that all
the relations between the generators of $\PB$ get mapped to relations
in $\PCD$ by $\Zc$. One can verify (using the Mac~Lane coherence
theorem~\cite{MacLane:Categories}) that the relations between the
generators of $\PB$ are the (repeated) $d_i$ images of the relations
(see also~\cite{Bar-Natan:NAT}):
\begin{itemize}
\item The pentagon $d_4a\circ d_2a\circ d_0a = d_1a\circ d_3a$,
  as above.
\item The hexagons $d_1\sigma^{\pm 1}=a\circ d_0\sigma^{\pm 1}\circ
  a^{-1}\circ d_3\sigma^{\pm 1}\circ a$, as above.
\item Locality in space: (slashes
  (\raisebox{0.5mm}{$\!\setlength{\unitlength}{0.5\standardunitlength}
	\begin{array}{c}  \hspace{-1.7mm}
        	\raisebox{-2pt}{\input draws/slash.tex }
        	\hspace{-1.9mm}
	\end{array}
\!\!$}) indicate bundles
  of strands)
  \[ \printname{LocSpace}
	\setlength{\unitlength}{0.5\standardunitlength}
	\begin{array}{c}  \hspace{-1.7mm}
        	\raisebox{-8pt}{\input draws/LocSpace.tex }
        	\hspace{-1.9mm}
	\end{array}
. \]
  Here $A$ and $B$ can each be either $a^{\pm 1}$ or $\sigma^{\pm 1}$.
\item Locality in scale:
  \[ \printname{LocRight}
	\setlength{\unitlength}{0.5\standardunitlength}
	\begin{array}{c}  \hspace{-1.7mm}
        	\raisebox{-8pt}{\input draws/LocRight.tex }
        	\hspace{-1.9mm}
	\end{array}
\quad\text{and}\quad\printname{LocSigma}
	\setlength{\unitlength}{0.5\standardunitlength}
	\begin{array}{c}  \hspace{-1.7mm}
        	\raisebox{-8pt}{\input draws/LocSigma.tex }
        	\hspace{-1.9mm}
	\end{array}
. \]
  Here $A$, $B$ and $C$ can each be either $a^{\pm 1}$ or $\sigma^{\pm 1}$.
\end{itemize}

Clearly, $\Zc$ respects the pentagon and the hexagons because $\Phi$
satisfies the pentagon and the hexagon axioms in the definition of
an associator. By Lemma~\ref{Locality}, $\Zc$ respects the locality
relations.  Hence $\Zc$ is well defined on morphisms of $\PB$. One can
verify that $\Zc({\bold I})\subset\F{1}\PCDc$, and hence $\Zc$ makes
sense on $\PBc$.  Finally, the fact that $\Zc$ intertwines the coproduct
$\Box$ and the operations $s_i$ follows from the group-like property
and the non-degeneracy of $\Phi$, respectively.

The proof in the case of $B^{(m)}$, $C^{(m)}$ and $\ASSm$ is essentially
identical.
\end{pf*}

\begin{proposition} \lbl{Invertible} Every structure preserving functor
$\Zc:\hat{B}\to\hat{C}$ or $\Zm:B^{(m)}\to C^{(m)}$ is invertible.
\end{proposition}

\begin{pf}
The unipotent completion $\widehat{PB}_n$ of the pure braid
group $PB_n$ on $n$ strands can be identified with the
the ring of morphisms in $\PBc$ from the $n$-point object
$O_r=(\bullet(\bullet\ldots(\bullet\bullet)\ldots))$ back to
itself that cover the identity permutation in $\PP$. Similarly,
$\Apbc_n$ can be identified with the ring of self-morphisms of $O_r$
in $\PCDc$ that cover the identity permutation. Thus, a functor
$\Zc:\widehat{B}\to\widehat{C}$ induces a filtration-preserving ring morphism
$\hat{Z}_n:\widehat{PB}_n\to\Apbc_n$.  In $\PBc$ ($\PCDc$) every morphism
can be written as a composition of invertible morphisms and an element
of $\widehat{PB}_n$ ($\Apbc_n$), and hence it is enough to prove that
$\hat{Z}_n$ is an isomorphism for every $n$. Finally, it is enough to
do that on the level of associated graded spaces. That is, we only need
to show that $\bar{Z}^{m}_n:\G{m}PB_n=I^m/I^{m+1}\to\G{m}\Apb_n$ is an
isomorphism for every $n$ and $m$, where $I$ is the augmentation ideal
of $PB_n$.

Let $\sigma^{ij}$ with $1\leq i<j\leq n$ be the standard generators of
$PB_n$:
\[ \sigma^{ij}=\printname{sigmaij}
	\setlength{\unitlength}{0.5\standardunitlength}
	\begin{array}{c}  \hspace{-1.7mm}
        	\raisebox{-8pt}{\input draws/sigmaij.tex }
        	\hspace{-1.9mm}
	\end{array}
. \]
In $\PB$, the parenthesized braid corresponding to $\sigma^{ij}$ is a
conjugate of an extension of $\sigma^2$:
\[ \printname{Conjugate}
	\setlength{\unitlength}{0.5\standardunitlength}
	\begin{array}{c}  \hspace{-1.7mm}
        	\raisebox{-8pt}{\input draws/Conjugate.tex }
        	\hspace{-1.9mm}
	\end{array}
. \]
As $\Zc(\sigma^2)=\tilde{R}^2=\exp H$ and $\Zc$ preserves all structure,
we find that
\[ \Zc(\sigma^{ij}-1)=\Zc(C)^{-1}(\exp t^{ij}-1)\Zc(C)=
  t^{ij}+(\text{higher terms}).
\]
Ergo,
\begin{equation} \lbl{Onto}
  \Zc\left(
    (\sigma^{i_1j_1}-1)(\sigma^{i_1j_2}-1)\cdots(\sigma^{i_kj_k}-1)
  \right)=
  t^{i_1j_1}t^{i_1j_2}\cdots t^{i_kj_k}+(\text{higher terms}),
\end{equation}
and the maps $\bar{Z}^m_n$ are surjective. Furthermore, as we mod out by
$I^{m+1}$, products of the form
\[ (\sigma^{i_1j_1}-1)(\sigma^{i_1j_2}-1)\cdots(\sigma^{i_mj_m}-1) \]
generate $\G{m}PB_n$, and hence it is enough to verify the injectivity of
$\bar{Z}^m_n$ on such products.

To do this we attempt to construct an inverse map $Y^m_n$ by setting
\[ Y^m_n(t^{i_1j_1}t^{i_1j_2}\cdots t^{i_mj_m})
  =(\sigma^{i_1j_1}-1)(\sigma^{i_1j_2}-1)\cdots(\sigma^{i_mj_m}-1).
\]
We only need to show that $Y^m_n$ is well defined; i.e., that it
carries the relations in $\G{m}\Apb_n$ to relations in $\G{m}PB_n$. This
is a routine verification. For example, if $i<j<k$, the braid relation
$[\sigma^{jk},\sigma^{ij}\sigma^{ik}]=0$ (``the third Reidemeister
move'') implies $[\sigma^{jk}-1,(\sigma^{ij}-1)+(\sigma^{ik}-1) -
(\sigma^{ij}-1)(\sigma^{ik}-1)]=0$. Notice that the last term in this
equality lies in a higher power of the augmentation ideal, and hence it
can be ignored. What remains proves that $Y^m_n$ maps the $4T$ relation
$[t^{jk},t^{ij}+t^{ik}]$ to 0 in the case when $i<j<k$.
\end{pf}

\begin{remark} In the language of Vassiliev
invariants, the last proof is essentially the identification of the space
of weight systems for pure braids with the dual of $\Apb$. If you know that
language, you may find it amusing to translate the above proof to the
Vassiliev setting.
\end{remark}

\begin{remark} Implicitly in the proof of Proposition~\ref{Invertible} we
have also proved that $\hat{C}$ is the ``associated graded mathematical
structure'' of the filtered structure $\hat{B}$.
\end{remark}

Propositions~\ref{KeyObservation} and~\ref{Invertible} imply the following:

\begin{theorem} The set $\ASSc$ ($\ASSm$) can be identified with the set
of all structure-preserving isomorphisms $\hat{B}\to\hat{C}$ ($B^{(m)}\to
C^{(m)}$). \qed
\end{theorem}

This would not be of much use if it was not for the following theorem, proven
by Drinfel'd~\cite{Drinfeld:QuasiHopf, Drinfeld:GalQQ} using
complex-analytic techniques:

\begin{theorem} \lbl{NonEmpty} The set $\ASSc({\Bbb C})$ (and thus $\ASSm$)
is non-empty. \qed
\end{theorem} 

This, in turn, allows us to use Principle~\ref{BasicPrinciple} and get:

\begin{theorem} \lbl{BasicTheorem} The groups $\GTc({\Bbb C})$ and
$\GRTc({\Bbb C})$ act simply transitively on $\ASSc({\Bbb C})$ on the right
and on the left respectively, and their actions commute. The same holds
for $\GTm({\Bbb C})$, $\GRTm({\Bbb C})$, and $\ASSm({\Bbb C})$. \qed
\end{theorem}

It is a consequence (and indeed, the purpose) of our main theorem
below, that Theorems~\ref{NonEmpty} and~\ref{BasicTheorem} also hold
over ${\Bbb Q}$.

%% file: draws/slash.tex
%
\begingroup\makeatletter\ifx\SetFigFont\undefined%
\gdef\SetFigFont#1#2#3#4#5{%
  \reset@font\fontsize{#1}{#2pt}%
  \fontfamily{#3}\fontseries{#4}\fontshape{#5}%
  \selectfont}%
\fi\endgroup%
\begin{picture}(174,189)(0,-10)
\thicklines
\path(87,12)(87,162)
\path(12,12)(162,162)
\end{picture}

%% file: theorem.tex
\section{The Main Theorem} \lbl{Main}

\subsection{The statement, consequences, and first reduction}

Our main theorem is:

\begin{theorem} \lbl{MainTheorem} (Proof on page~\pageref{MainProof})
The natural map $\ASSm({\Bbb C})\to\ASSmm({\Bbb C})$ is surjective.
\end{theorem}

This theorem means that an associator can be constructed degree by
degree. Furthermore, if $\Phi_{m-1}\in\ASSmm$ is an associator up to
degree $m-1$ and $\Phi_m=\Phi_{m-1}+\varphi_m$, with $\deg\varphi_m=m$,
then the equations\footnote{More on these equations can be found in
Drinfel'd~\cite{Drinfeld:GalQQ} and in~\cite{Bar-Natan:NAT}.} that
$\varphi_m$ has to satisfy for $\Phi_m$ to be an associator up to degree
$m$ are non-homogeneous linear, with a constant term determined
algebraically from $\Phi_{m-1}$. Therefore, if a $\Phi_{m-1}$ is found
over the rationals, then a $\Phi_m$ can be found over the rationals
(i.e., the statement of Theorem~\ref{MainTheorem} also holds over ${\Bbb
Q}$). Proceeding using induction, we find that a rational associator
exists (and so Theorems~\ref{NonEmpty} and~\ref{BasicTheorem} also hold
over ${\Bbb Q}$).

\begin{corollary} Rational associators exist and can be constructed
iteratively. \qed
\end{corollary}

Let $P$ be the automorphism of $\Apb$ that sends every generator $t^{ij}$
to its negative $-t^{ij}$. It is clear that $P$ preserves $\ASSm$ (it
simply switches the positive and negative hexagon identities while not
touching the pentagon identity). If $\Phi_{m-1}\in\ASSmm$ is even (i.e.,
satisfies $\Phi_{m-1}=P\Phi_{m-1}$), it can be lifted to an even
$\Phi_m\in\ASSm$: Simply take any lifting $\Phi'_m$ and set
$\Phi_m=(\Phi'_m+P\Phi'_m)/2$. This is an associator because the set of
liftings of $\Phi_{m-1}$ is affine, as it is determined by the solutions of
a non-homogeneous linear equation.

\begin{corollary} Rational even associators exist and can be constructed
iteratively. \qed
\end{corollary}

\begin{remark} Even associators were given a topological interpretation
in~\cite{LeMurakami:Parallel}, and were used further
in~\cite{LeMurakamiOhtsuki:Universal}.
\end{remark}

\begin{lemma} \lbl{Enough}
To prove Theorem~\ref{MainTheorem} it is enough to prove that the natural
homomorphism $\GRTm({\Bbb C})\to\GRTmm({\Bbb C})$ is surjective.
\end{lemma}

\begin{pf} By Theorem~\ref{NonEmpty}, $\ASSm({\Bbb C})$ is non-empty,
and so there exists at least one $\Phi_{m-1}\in\ASSmm({\Bbb C})$ that
extends to a $\Phi_m\in\ASSm({\Bbb C})$. Take now any other element
$\Phi'_{m-1}$ of $\ASSmm({\Bbb C})$. By Theorem~\ref{BasicTheorem}, it
can be pushed to $\Phi_{m-1}$ by some $G_{m-1}\in\GRTmm({\Bbb C})$. Take
a $G_m\in\GRTm({\Bbb C})$ that extends $G_{m-1}$, and use it to pull
$\Phi_m$ back to become an extension $G_m^{-1}\Phi_m$ of $\Phi'_{m-1}$,
as required.
\end{pf}

\subsection{More on the group $\protect\GRTc$}

To prove the surjectivity of $\GRTm(A)\to\GRTmm(A)$ for some ground
algebra $A$, we need to know some more about $\GRTm=\Aut C^{(m)}$ and
about the structure $C^{(m)}$ itself. Recall that the category $\PCD$
is generated by the (repeated) $d_i$ images of the specific morphisms
$a^{\pm 1}$, $X$ and $H$.

\begin{proposition} \lbl{PaCDRels} The (repeated) $d_i$ images of the
relations below generate all the relations between generators of $\PCD$:
\begin{itemize}
\item $X$ is its own inverse and it commutes with $H$.
\item The pentagon $d_4a\circ d_2a\circ d_0a = d_1a\circ d_3a$,
  as for the category $\PB$.
\item The classical hexagon
  \begin{equation} \lbl{ClassicalHexagon}
    \printname{ClassicalHexagon}
	\setlength{\unitlength}{0.5\standardunitlength}
	\begin{array}{c}  \hspace{-1.7mm}
        	\raisebox{-8pt}{\input draws/ClassicalHexagon.tex }
        	\hspace{-1.9mm}
	\end{array}
\quad;\qquad
    d_1X = a\circ d_0X\circ a^{-1}\circ d_3X\circ a.
  \end{equation}
\item The semi-classical hexagon (the name is explained in
  Remark~\ref{QuantumHexagon})
  \begin{equation} \lbl{SCHexagon}
    d_1\left(\,\printname{HX}
	\setlength{\unitlength}{0.35\standardunitlength}
	\begin{array}{c}  \hspace{-1.7mm}
        	\raisebox{-8pt}{\input draws/HX.tex }
        	\hspace{-1.9mm}
	\end{array}
\right)\eqdef\printname{SCHexLHS}
	\setlength{\unitlength}{0.5\standardunitlength}
	\begin{array}{c}  \hspace{-1.7mm}
        	\raisebox{-8pt}{\input draws/SCHexLHS.tex }
        	\hspace{-1.9mm}
	\end{array}

    =\printname{SCHexRHS}
	\setlength{\unitlength}{0.5\standardunitlength}
	\begin{array}{c}  \hspace{-1.7mm}
        	\raisebox{-8pt}{\input draws/SCHexRHS.tex }
        	\hspace{-1.9mm}
	\end{array}
;
  \end{equation}
\[ d_1H\circ d_1X=
  a\circ d_0H \circ d_0X \circ a^{-1}\circ d_3X \circ a
  + a\circ d_0X\circ a^{-1}\circ d_3H \circ d_3X \circ a.
\]
\item Locality in space as in $\PB$ (but with $A,B\in\{a^{\pm 1},X,H\}$).
\item Locality in scale:
  \[ \printname{LocRight}
	\setlength{\unitlength}{0.5\standardunitlength}
	\begin{array}{c}  \hspace{-1.7mm}
        	\raisebox{-8pt}{\input draws/LocRight.tex }
        	\hspace{-1.9mm}
	\end{array}
,\qquad\qquad\printname{LocX}
	\setlength{\unitlength}{0.5\standardunitlength}
	\begin{array}{c}  \hspace{-1.7mm}
        	\raisebox{-8pt}{\input draws/LocX.tex }
        	\hspace{-1.9mm}
	\end{array}
, \]
  \[ \text{and}\qquad\qquad\printname{LocH}
	\setlength{\unitlength}{0.5\standardunitlength}
	\begin{array}{c}  \hspace{-1.7mm}
        	\raisebox{-8pt}{\input draws/LocH.tex }
        	\hspace{-1.9mm}
	\end{array}
, \]
  with $A,B,C\in\{a^{\pm 1},X,H\}$.
\end{itemize}
\end{proposition}

\begin{pf} (sketch)
Let $\TMP$ be the fibered-linear category freely generated by the repeated
$d_i$ images of $a^{\pm 1}$, $X$ and $H$ in $\PCD$, modulo the relations
listed above. There is an obvious functor ${\bold F}:\TMP\to\PCD$, which
is well defined because the relations above are indeed relations in $\PCD$
($4T$ is needed to verify the third locality in scale relation with $A$
or $B$ equal $H$). The category $\TMP$ is graded by declaring that $\deg
a=\deg X=0$, that $\deg H=1$, and that the operations $d_i$ preserve
degree. Clearly, the functor $\bold F$ preserves degrees.  We need to
show that ${\bold F}$ is invertible, and we do so by constructing an
inverse ${\bold G}:\PCD\to\TMP$ in steps as follows:

\begin{myenumi}
\def\theenumi{(\arabic{enumi})} \setcounter{enumi}{1}
\def\theenum{\theenumi\refstepcounter{enumi}}
\item[\theenum]
  There is no problem with constructing ${\bold G}$ in degree
  $0$.  The relevant generators of $\TMP$ are commutativities $X$ and
  associativities $a^{\pm 1}$, and the relevant relations are (some of)
  the locality relations and the pentagon and the classical hexagon. Thus
  the existence of ${\bold G}$ in degree $0$ is exactly the Mac~Lane
  coherence Theorem~\cite{MacLane:Categories}.
\end{myenumi}

\par\noindent
Let $\PCD_r$ ($\TMP_r$) be the algebra of self-morphisms of the
object $O_r=(\bullet(\bullet\ldots(\bullet\bullet)\ldots))$ in $\PCD$
($\TMP$) that cover the identity permutation in $\PP$, and let ${\bold
F}_r:\TMP_r\to\PCD_r$ be the obvious restriction of ${\bold F}$. Our next
objective is to construct ${\bold G}_r$, an inverse of ${\bold
F}_r$. There is no loss of generality in assuming that all morphisms that
we deal with involve exactly $n$ strands (for some fixed $n$). With this
in mind, the $\PCD_r$ can be identified with $\Apb_n$.

\begin{myenumi}
\def\theenumi{(\arabic{enumi})} \setcounter{enumi}{2}
\def\theenum{\theenumi\refstepcounter{enumi}}
\item[\theenum]
  Construct ${\bold G}_r$ in degree $1$. It is enough to specify the
  image in $\TMP_r$ of $t^{ij}\in\Apb_n$, and to check that ${\bold
  G}_r$ is indeed the inverse of ${\bold F}_r$ in degree $1$. So for
  $i<j$ set ${\bold G}(t^{ij})=P_{ij}^{-1}\circ H_n\circ P_{ij}$. Here
  $H_n=d^{n-1}_0d^{n-2}_0\cdots d^2_0H$ is $H$ extended by adding strands
  on the left and $P_{ij}\in\mor_{\TMP}(O_r,O_r)$ with $\deg P_{ij}=0$
  corresponds to the parenthesized permutation that takes the $j$th
  strand to be the last and the $i$th to be the one before the last,
  while preserving the order of all other strands. (Step (1) implies
  that it does not matter which particular generator combination
  we choose for $P_{ij}$).  Now ${\bold F}_r\circ{\bold G}_r=Id$ is
  trivial, and ${\bold G}_r\circ{\bold F}_r=Id:\TMP_r\to\TMP_r$ is not
  hard to check.  Indeed, a degree $1$ morphism in $\TMP_r$ contains
  exactly one (repeated $d_i$ image of) $H$. By the semi-classical
  hexagon we can replace cabled $H$'s by extended ones (terminology
  as in Definition~\ref{PCDOperators}), and extended $H$'s can be slid
  right using locality in scale relations:
  \[ \printname{SlideRight}
	\setlength{\unitlength}{0.5\standardunitlength}
	\begin{array}{c}  \hspace{-1.7mm}
        	\raisebox{-8pt}{\input draws/SlideRight.tex }
        	\hspace{-1.9mm}
	\end{array}
. \]
  Finally, on ``right justified'' $H$'s there is almost nothing to prove.
\item[\theenum]
  By extending ${\bold G}_r$ multiplicatively to higher degrees, we find
  that the free algebra $FT$ generated by $\G{1}\TMP_r$ is isomorphic
  to the free algebra $FA$ generated by $\G{1}\Apb_n$.  The algebra
  $\TMP_r$ is the quotient of $FT$ by the quadratic locality relations:
  locality in space with $A=B=H$, and the third locality in scale relation
  with either $A=H$ or $B=H$.  The algebra $\G{1}\Apb_n$ is the quotient
  of $FA$ by quadratic relations: the relation $[t^{ij},t^{kl}]=0$ and
  the $4T$ relation $[t^{jk},t^{ij}+t^{ik}]=0$.  Quite clearly, these
  relations correspond under the isomorphisms between $FT$ and $FA$;
  the locality relation $\printname{Suggestive4T}
	\setlength{\unitlength}{0.5\standardunitlength}
	\begin{array}{c}  \hspace{-1.7mm}
        	\raisebox{-8pt}{\input draws/Suggestive4T.tex }
        	\hspace{-1.9mm}
	\end{array}
$, for example, is
  sent to the $4T$ relation. We conclude that the quotients $\TMP_r$
  and $\Apb_n$ are isomorphic via ${\bold F}_r$ and ${\bold G}_r$.
\end{myenumi}

\par\noindent
Finally, we get back to constructing ${\bold G}$:

\begin{myenumi}
\def\theenumi{(\arabic{enumi})} \setcounter{enumi}{4}
\def\theenum{\theenumi\refstepcounter{enumi}}

\item[\theenum]
  Every morphism $M$
  in $\PCD$ can be written uniquely as a composition $P_1\circ D\circ
  P_2$ where $D\in\Apb_n=\PCD_r$, $P_{1,2}$ are of degree $0$, and
  $P_1$ induces the identity permutation (between possibly different
  parenthesizations). Define ${\bold G}(M)={\bold G}(P_1)\circ
  {\bold G}_r(D)\circ {\bold G}(P_2)$. Clearly, ${\bold G}$ is the
  inverse of ${\bold F}$.
  \qed
\end{myenumi}
\renewcommand{\qed}{}
\end{pf}

\begin{remark} \lbl{QuantumHexagon}
Let $\epsilon$ be a formal parameter satisfying $\epsilon^2=0$, and
let $\PCD_\epsilon$ be defined as $\PCD$, only with coefficients in
the algebra $A[\epsilon]$ rather than the algebra $A$. Let $R_\epsilon$
be the morphism $(\exp\epsilon H)\circ X$ in $\PCD_\epsilon$, and consider
the ``quantum'' hexagon relation for $R_\epsilon$:
\[ d_1R_\epsilon
    = a\circ d_0R_\epsilon\circ a^{-1}\circ d_3R_\epsilon\circ a.
\]
A quick visual inspection of equations~\eqref{HexagonPlus} (with
$R_\epsilon$ replacing $\sigma$),~\eqref{ClassicalHexagon}
and~\eqref{SCHexagon} reveals that the classical and semi-classical
hexagon relations are the degree $0$ and $1$ parts (in $\epsilon$) of
the quantum hexagon relation, explaining their names.
\end{remark}

\begin{remark} \lbl{CablingRelation}
Modulo the other relations, the semi-classical hexagon is equivalent to the
simpler but less conceptual ``cabling relation'', $d_2H = a^{-1}\circ
d_3H\circ a + d_0X\circ a^{-1}\circ d_3H\circ a\circ d_0X$:
\[ d_2\left(\,\printname{H}
	\setlength{\unitlength}{0.5\standardunitlength}
	\begin{array}{c}  \hspace{-1.7mm}
        	\raisebox{-8pt}{\input draws/H.tex }
        	\hspace{-1.9mm}
	\end{array}
\right)\eqdef\printname{CablingLHS}
	\setlength{\unitlength}{0.5\standardunitlength}
	\begin{array}{c}  \hspace{-1.7mm}
        	\raisebox{-8pt}{\input draws/CablingLHS.tex }
        	\hspace{-1.9mm}
	\end{array}

  =\printname{CablingRHS}
	\setlength{\unitlength}{0.5\standardunitlength}
	\begin{array}{c}  \hspace{-1.7mm}
        	\raisebox{-8pt}{\input draws/CablingRHS.tex }
        	\hspace{-1.9mm}
	\end{array}
.
\]
\end{remark}

By Claim~\ref{PaCDGens} and Remark~\ref{FixXH},
any $G\in\GRTm$ is determined by its action on the generator $a$ of
$\PCDm$, and thus it is determined by the unique $\Gamma\in\Apbm_3$ for which
$G(a)=\Gamma\cdot a$. Just as in the proof of
Proposition~\ref{KeyObservation}, the relations of
Proposition~\ref{PaCDRels} impose relations on $\Gamma$:

\begin{proposition} \lbl{GammaEquations}
The group $\GRTm$ ($\GRTc$) can be identified (as a set) with the set
of all group-like non-degenerate $\Gamma\in\Apbm_3$ ($\Gamma\in\Apbc_3$)
satisfying:
\begin{itemize}
\item The pentagon equation $d_4\Gamma\cdot d_2\Gamma\cdot d_0\Gamma =
  d_1\Gamma\cdot d_3\Gamma$.
\item The classical hexagon equation
  $1=\Gamma\cdot(\Gamma^{-1})^{132}\cdot\Gamma^{312}$.
\item The semi-classical hexagon equation
  \[ d_1t^{12}=\Gamma\cdot
    \left(t^{23}\cdot(\Gamma^{-1})^{132}+(\Gamma^{-1})^{132}\cdot t^{13}\right)
    \cdot\Gamma^{312},
  \]
  or, equivalently, the cabling equation $d_2t^{12}=\Gamma^{-1}\cdot
  t^{12}\cdot\Gamma + (\Gamma^{-1}\cdot t^{12}\cdot\Gamma)^{132}$.
\end{itemize}
\end{proposition}

\begin{pf} The group-like property and the non-degeneracy of $\Gamma$
correspond to the fact that $G$ preserves $\Box$ and the operations
$s_i$. The pentagon, classical and semi-classical hexagon, and cabling
equations correspond to their namesakes in Proposition~\ref{PaCDRels}. The
other relations in Proposition~\ref{PaCDRels} impose no further constrains
on $\Gamma$; the locality relations follow from Lemma~\ref{Locality}
and the relations $X^2=1$ and $XH=HX$ do not involve $\Gamma$ at all.
\end{pf}

\begin{warning} The product of $\GRTm$ ($\GRTc$) is not the product of
$\Apbm$ ($\Apbc$). See Proposition~\ref{GRTGroupLaw}.
\end{warning}

\begin{remark} The classical hexagon axiom for $\Gamma\in\GRTm$ implies
that $\Gamma=1+$(higher degree terms).
\end{remark}

\begin{remark} \lbl{Spirit}
In the spirit of Remark~\ref{QuantumHexagon}, the classical and
semi-classical hexagon equations can be replaced by a single ``quantum
hexagon equation'' written in $\Apbm_3(A[\epsilon])$:
\begin{equation} \lbl{QuantumHexagonEq}
  e^{\epsilon(t^{13}+t^{23})} =
  \Gamma\cdot e^{\epsilon t^{23}}\cdot(\Gamma^{-1})^{132}\cdot
    e^{\epsilon t^{13}}\cdot\Gamma^{312}.
\end{equation}
\end{remark}

\subsection{The second reduction}

\begin{theorem} \lbl{SCFollows} (Proof on page~\pageref{SCFollowsProof})
The pentagon and classical hexagon equations for $\Gamma\in\Apbm_3$ imply
the semi-classical hexagon equation (and hence the cabling equation).
\end{theorem}

Assuming Theorem~\ref{SCFollows}, the proof of Theorem~\ref{MainTheorem}
reduces to an easy observation and some standard (but non-trivial) facts
from the theory of affine group schemes.

\begin{pf*}{Proof of Theorem~\ref{MainTheorem}} \lbl{MainProof}
By Lemma~\ref{Enough}, it is enough to show that the natural homomorphism
$\pi:\GRTm({\Bbb C})\to\GRTmm({\Bbb C})$ is surjective. In the next
paragraph we will show that $\pi$ is a homomorphism of connected reduced
algebraic group schemes. Hence it is enough to prove this statement
at the level of Lie algebras, and the Lie algebras are given by the
linearizations near the identity $1$ of the defining equations, the
pentagon and the classical hexagon. These linearizations are
\begin{equation} \lbl{grt}
  d_4\gamma + d_2\gamma + d_0\gamma = d_1\gamma + d_3\gamma
  \qquad\text{and}\qquad
  0=\gamma - \gamma^{132} + \gamma^{312}.
\end{equation}
Clearly, any solution to degree $m-1$ of these equations can be extended
to a solution to degree $m$ (for example, by taking the degree $m$ piece
to be $0$). Notice that if the cabling relation was still present, this
would not have been so easy: The linearization of the cabling relation is
$0=[t^{12},\gamma]+[t^{13},\gamma^{132}]$, and this equation at degree
$m$ imposes a (possibly new) condition on the degree $m-1$ piece of
$\gamma$.

All that is left now is some standard algebraic geometry. We defined
$\GRTm(A)$ for an arbitrary ground algebra $A$ in a functorial
way, and saw that it is always defined by the same equations
(Proposition~\ref{GammaEquations}). Thus $\GRTm$ (regarded as a functor
from the category of ${\Bbb Q}$-algebras to the category of groups) is an
affine group scheme (see e.g.~\cite[section~1.2]{Waterhouse:GroupSchemes})
for any $m$ (and similarly, the map $\GRTm\to\GRTmm$ is a homomorphism
of affine group schemes). $\GRTm$ has a faithful representation
in the vector space $V$ of parenthesized chord diagrams whose
skeleton is $a$ (already the action of $G\in\GRTm$ on $a$ itself
determines $G$). Thus $\GRTm$ can be regarded as an algebraic
matrix group. Notice that for any $G\in\GRTm$, we have $G(X)=X$,
$G(H)=H$, and $G(a)=a+(\text{higher degrees})$, and hence for any
homogeneous $v\in V$ we have $G(v)=v+(\text{higher degrees})$. Hence
$G$ is unipotent, and $\GRTm$ is a unipotent group~\cite[section
8]{Waterhouse:GroupSchemes}. As we are working in characteristic
$0$, $\GRTm$ is reduced~\cite[section 11.4]{Waterhouse:GroupSchemes}
(and hence~\eqref{grt} defines its Lie algebra) and $\GRTmm$ is
connected~\cite[section 8.5]{Waterhouse:GroupSchemes}.
\end{pf*}

\begin{remark} Very little additional effort as in the paragraph
following Theorem~\ref{MainTheorem} shows that $\GRTm(A)\to \GRTmm(A)$
is surjective for any $A$.
\end{remark}

\subsection{A cohomological interlude} Before we can prove
Theorem~\ref{SCFollows}, we need to know a bit about the second
cohomology of $\Apb_n$. There are two relevant ways of
turning the list $\Apb_2,\ \Apb_3,\ldots$ into a cochain
complex. The first is to define $d=d^n:\Apb_n\to\Apb_{n+1}$
by $d^n=\sum_{i=0}^{n+1}(-1)^id^n_i$. The second is to define
$\tilde{d}=\tilde{d}^n:\Apb_{n+1}\to\Apb_{n+2}$ (notice the shift
in dimension) by $\tilde{d}^n=\sum_{i=0}^{n+1}(-1)^i\tilde{d}^n_i$,
where $\tilde{d}_i=\tilde{d}^n_i=d^{n+1}_i$ for $i\leq n$, and
$\tilde{d}_{n+1}=\tilde{d}^n_{n+1}=(d^{n+1}_{n+2})^{12\ldots n(n+2)(n+1)}$
is the operation of ``adding an empty strand between strands $n$ and
$n+1$'':
\[ \tilde{d}_3\left(\printname{Psi}
	\setlength{\unitlength}{0.5\standardunitlength}
	\begin{array}{c}  \hspace{-1.7mm}
        	\raisebox{-8pt}{\input draws/Psi.tex }
        	\hspace{-1.9mm}
	\end{array}
\right)=\printname{Psi124}
	\setlength{\unitlength}{0.5\standardunitlength}
	\begin{array}{c}  \hspace{-1.7mm}
        	\raisebox{-8pt}{\input draws/Psi124.tex }
        	\hspace{-1.9mm}
	\end{array}
. \]

For the purpose of proving Theorem~\ref{SCFollows}, all we need is to
understand $H^2_{\tilde{d}}$:

\begin{proposition} \lbl{CohomResult} $H^2_{\tilde{d}}$ is 2-dimensional
and is generated by $t^{12}$ (in degree $1$) and $[t^{13},t^{23}]$ (in
degree 2).
\end{proposition}

\begin{pf}

It is well known (see e.g.~\cite{Kohno:MonRep, Drinfeld:GalQQ,
Hutchings:SingularBraids, Bar-Natan:Braids}) that as vector spaces,
$\Apb_{n+1}=\Apb_n\otimes {\cal T}V^n$, where ${\cal T}V^n$ denotes
the tensor algebra on the $n$-dimensional vector space $V^n$ generated
by $t^{1(n+1)},\ldots,t^{n(n+1)}$ (as algebras, this is a semi-direct
product).  Furthermore, $\tilde{d}^n_i$ and the strand removal operations
$\tilde{s}^n_i\eqdef s^{n+1}_i$ preserve this decomposition, and define
a structure of a cosimplicial vector space on each of $\Apb_{n+1}$,
$\Apb_n$, and $V^n$. The cosimplicial structure induced on $\Apb_n$
coincides with the one it already has $((d^n_i,s^n_i))$, and hence by
the Eilenberg-Zilber theorem and the K\"{u}nneth formula
\begin{equation} \lbl{decomposition}
  H^\star_{\tilde{d}}=H^\star_d\hat{\otimes}\hat{\cal T}H^\star(V^\star).
\end{equation}
(Here $\hat{\otimes}$ denotes the ${\Bbb Z}/2{\Bbb Z}$-graded tensor
product and $\hat{\cal T}$ denotes the tensor algebra formed using
$\hat{\otimes}$).

\par\noindent
{\em Computing $H^\star_d$:} The cohomology $H^\star_d$ is very
hard to compute. Indeed, if we could compute $H^4_d$, we probably
needn't have written this paper at all (see~\cite{Drinfeld:GalQQ,
Bar-Natan:Braids}). But up to $H^2_d$ there is no difficulty in computing
by hand. The algebras $\Apb_0$ and $\Apb_1$ contain only multiples of
the identity element. The algebra $\Apb_2$ contains only the powers
of $t^{12}$. The differential $d^0:\Apb_0\to\Apb_1$ is the zero map,
the differential $d^1:\Apb_1\to\Apb_2$ is injective, mapping the
identity of $\Apb_1$ to the identity of $\Apb_2$. Finally, let us study
$d^2(t^{12})^m\in\Apb_3$. Setting $c=t^{12}+t^{13}+t^{23}\in\Apb_3$,
we get:
\[ d^2(t^{12})^m
  = \sum_{i=0}^3 (-1)^id^2_i(t^{12})^m
  = (t^{23})^m-(c-t^{12})^m+(c-t^{23})^m+(t^{12})^m.
\]
The relations of Definition~\ref{def:Apb} (in the case $n=3$) can be
rewritten in terms of the new generators $t^{12}$, $t^{23}$ and $c$
of $\Apb_3$. In these terms, they are equivalent to the statement
``$c$ is central''. Thus $\Apb_3$ is the central extension by $c$ of
the free algebra in $t^{12}$ and $t^{23}$. Looking at the coefficient
of (say) $c(t^{12})^{(m-1)}$ in $d^2(t^{12})^m$ as computed above, we
find that $d^2(t^{12})^m\neq 0$ for $m\geq 2$. It is easy to verify that
$d^2(t^{12})^m=0$ for $m=0,1$. In summary, we found that $\dim H^0_d=1$,
with the generator being the unit of $\Apb_0$, that $\dim H^1_d=0$,
and that $H^2_d$ is one dimensional and is generated by $t^{12}$.

\par\noindent
{\em Computing $H^\star(V^\star)$:} By the normalization theorem for
simplicial cohomology the complex $(V^n)$ has the same cohomology as the
complex $(\hat{V}^n)$ defined by $\hat{C}^n=\bigcap_i\ker \tilde{s}^n_i$.
But it is clear that $\hat{V}^n=0$ unless $n=1$, and that $\tilde{C}^1$
is 1-dimensional. Thus $H^\star(V^\star)$ has only one generator,
$t^{12}$ in $H^1(V^\star)$. (The same computation appears in~\cite[Lemma
4.14]{Bar-Natan:NAT}).

\par\noindent
{\em Assembling the results:} Using~\eqref{decomposition} and the
above two cohomology computations, we find that $H^2_{\tilde{d}}$
is generated by the class of $t^{12}$ (coming from $H^2_d$) and a
degree 2 class coming from the class $t^{12}\otimes t^{12}$ in
$H^1(V^\star)\hat{\otimes}H^1(V^\star)$ via the K\"{u}nneth map. An
explicit computation of the latter (or a direct computation of the cycles
and boundaries, which is easy in this low dimension), shows that it is
the class of $[t^{13},t^{23}]$.
\end{pf}

\subsection{Proof of the semi-classical hexagon equation}

\begin{pf*}{Proof of Theorem~\ref{SCFollows}} \lbl{SCFollowsProof} {
\def\s{\clubsuit}
\def\G{\Gamma}
\def\h{\!\!\hexagon\!\!_\epsilon}
\def\p{{\!\!\pentagon\!\!}}
\def\e#1{e^{\epsilon t^{#1}}}
\def\em#1{e^{-\epsilon t^{#1}}}
\def\E#1{e^{\epsilon(#1)}}
\def\Em#1{e^{-\epsilon(#1)}}
\def\d#1{\tilde{d}_{#1}}
\def\D#1{d_{#1}}
\def\i#1{{(#1)^{-1}}}
Assume that for some $\Gamma\in\Apbm_3$ the pentagon and
the classical hexagon hold, but the semi-classical hexagon
doesn't. By Remark~\ref{Spirit}, we know that the quantum
hexagon~\eqref{QuantumHexagonEq} has an error proportional to $\epsilon$.
Let $\epsilon\psi'$ be that error:
\[ 1+\epsilon\psi' =
  \Gamma\cdot e^{\epsilon t^{23}}\cdot(\Gamma^{-1})^{132}\cdot
  e^{\epsilon t^{13}}\cdot\Gamma^{312}\cdot e^{-\epsilon(t^{13}+t^{23})}.
\]
By assumption, $\psi'\neq 0$. Let $\psi$ be the lowest degree piece of
$\psi'$, and let $k=\deg\psi$. Clearly, $k\geq 2$. From this point on,
mod out by degrees higher than $k$.

We claim that
\begin{equation} \lbl{P20eq}
  \tilde{d}^2\psi=0.
\end{equation}
The proof of~\eqref{P20eq} is essentially contained in
Figure~\ref{P20}. How polyhedra correspond to identities of this kind was
explained in~\cite{Drinfeld:QuasiHopf}, and again in~\cite{Bar-Natan:NAT},
where the very same polyhedron appeared in a very similar context. For
completeness, we include the explanation here, in a very concrete form.
In Figure~\ref{P20} every edge is oriented and is labeled by some
invertible element of $\Apbm_4(A[\epsilon])$. There are 12 faces in
the figure (including the face at infinity). Each one corresponds to a
certain product in $\Apbm_4(A[\epsilon])$ by starting at the $\clubsuit$
symbol, going counterclockwise,  and multiplying the elements seen on the
edges (or their inverses depending on the edge orientations). These
products turn out to all be locality relations, or pentagons, or quantum
hexagons (or a permutation or a cabling/extension operation applied to
a pentagon or a quantum hexagon), as marked within each face.

\begin{figure}[htpb]
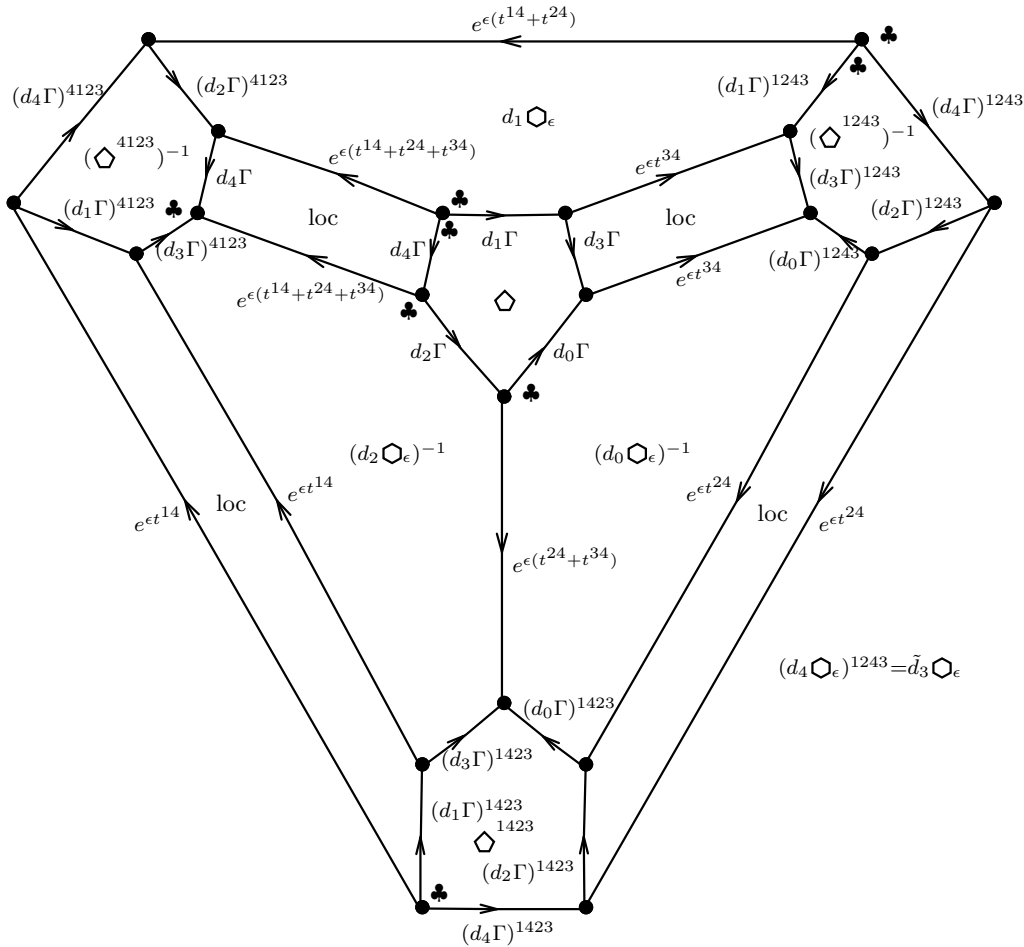

\[ \ \qquad\printname{P20}
	\setlength{\unitlength}{0.85\standardunitlength}
	\begin{array}{c}  \hspace{-1.7mm}
        	\raisebox{-8pt}{\input draws/P20.tex }
        	\hspace{-1.9mm}
	\end{array}
 \]
\caption{The proof of equation~\eqref{P20eq}.}
\lbl{P20} \end{figure}

For example (remember that we are ignoring degrees higher than $k$),
\[ \begin{array}{ccrrcl}
  \printname{p20-1}
	\setlength{\unitlength}{0.1\standardunitlength}
	\begin{array}{c}  \hspace{-1.7mm}
        	\raisebox{-8pt}{\input draws/p20-1.tex }
        	\hspace{-1.9mm}
	\end{array}
 & \hspace{-4mm}\rightarrow\!\!\! & \p: &
    1 & \!\!\!=\!\!\! & \D4\G\D2\G\D0\G\i{\D3\G}\i{\D1\G}, \\
  &&&&& \\
  \printname{p20-3}
	\setlength{\unitlength}{0.1\standardunitlength}
	\begin{array}{c}  \hspace{-1.7mm}
        	\raisebox{-8pt}{\input draws/p20-3.tex }
        	\hspace{-1.9mm}
	\end{array}
 & \hspace{-4mm}\rightarrow\!\!\! & \D1\h: &
    1\!+\!\D1\psi & \!\!\!=\!\!\! & \D1\G\e{34}\i{(\D1\G)^{1243}}
      \E{t^{14}+t^{24}}(\D2\G)^{4123}\Em{t^{14}+t^{24}+t^{34}}, \\
  &&&&& \\
  \printname{p20-2}
	\setlength{\unitlength}{0.1\standardunitlength}
	\begin{array}{c}  \hspace{-1.7mm}
        	\raisebox{-8pt}{\input draws/p20-2.tex }
        	\hspace{-1.9mm}
	\end{array}
 & \hspace{-4mm}\rightarrow\!\!\! & \D2\h^{-1}: &
    1\!-\!\D2\psi & \!\!\!=\!\!\! & (\text{product around shaded area}).
\end{array} \]

Combining these equations along the common edges we get
\[ \printname{p20-4}
	\setlength{\unitlength}{0.1\standardunitlength}
	\begin{array}{c}  \hspace{-1.7mm}
        	\raisebox{-8pt}{\input draws/p20-4.tex }
        	\hspace{-1.9mm}
	\end{array}
 \rightarrow 1+\D1\psi-\D2\psi
  =(\text{product around shaded area}).
\]
Continuing along the same line, we find that the product around the
whole figure is $1-\D0\psi+\D1\psi-\D2\psi$. On the other hand, this
product is itself a variant of the quantum hexagon --- $\i{\d3\h}$, as
marked on the face at infinity. So we learn that
$1-\D0\psi+\D1\psi-\D2\psi=1-\d3\psi$. But this is exactly~\eqref{P20eq}.

By~\eqref{P20eq} and Proposition~\ref{CohomResult}, we see that
if $k>2$ then $\psi$ must be in $\tilde{d}^1{\cal G}_k\Apb_2$. That
is, it must be a multiple of $\chi=\tilde{d}^1(t^{12})^k$. But as
$\Gamma$ is group-like, $\psi$ must be primitive: $\Box\psi=\psi\otimes
1+1\otimes\psi$. One easily verifies that $\chi$ is not primitive,
and hence $\psi=0$ as required. If $k=2$, equation~\eqref{P20eq}
and Proposition~\ref{CohomResult} tell us that $\psi$ is of the form
$c_1\tilde{d}^1(t^{12})^2+c_2[t^{13},t^{23}]$. A routine verification
shows that if the semi-classical hexagon relation is pre-multiplied
by $d_3 X$ and post-multiplied by $d_0 X$, then modulo the other
relations, it does not change. This means that $\psi^{213}=\psi$ (this
identity follows more easily from the cabling relation), and thus
$c_2=0$. But then the primitivity of $\psi$ implies that $c_1$ vanishes as
well, and thus $\psi=0$ as required.
}\end{pf*}

%% file: Just.tex
\section{Just for completeness} \lbl{Just}

For completeness, this section contains a description of the group
law of $\GRTc$, a description of its action on $\ASSc$, and similar
descriptions for the group $\GTc$. This information is not needed in the
main part of this paper. Throughout this section one can replace unipotent
completions by unipotent quotients ($\GRTm$, $\ASSm$, $\Apbm$, etc.)
with no change to the results.

\begin{proposition} \lbl{GRTGroupLaw}
The group law $\times$ of $\GRTc$ is expressed in terms of the $\Gamma$'s
(of Proposition~\ref{GammaEquations}) as
\begin{equation} \lbl{GammaProduct}
  \Gamma_1\times\Gamma_2=\Gamma_1\cdot\left(\left.\Gamma_2\right|_{
    t^{12}\to\Gamma_1^{-1}t^{12}\Gamma_1,\
    t^{13}\to(\Gamma_1^{-1})^{132}t^{13}\Gamma_1^{132},\
    t^{23}\to t^{23}
  }\right),
\end{equation}
where ``$\cdot$'' is the product of $\Apbc$, $\Gamma_1^{-1}$ is
interpreted in $\Apbc$, and the substitution above means: replace every
occurrence of $t^{12}$ in $\Gamma_2$ by $\Gamma_1^{-1}t^{12}\Gamma_1$,
etc.\ (In particular, we claim that this substitution is well defined
on $\Apbc$).
\end{proposition}

\begin{pf} $\Apbc_3$ can be identified with the algebra of
self-morphisms in $\PCDc$ of the object $(\bullet(\bullet\bullet))$. Let
$\overline\Gamma$ denote the self-morphism corresponding to a
$\Gamma\in\Apbc_3$. We have $\Gamma\cdot a=a\circ\overline\Gamma$, and
hence (with $\Gamma\mapsto G_\Gamma$ denoting the identification in
Proposition~\ref{GammaEquations})
\begin{equation} \lbl{GammaComposition}
  a\!\circ\!\overline{\Gamma_1\!\times\!\Gamma_2}
  = G_{\Gamma_1\!\times\!\Gamma_2}(a)
  = G_{\Gamma_1}(G_{\Gamma_2}(a))
  = G_{\Gamma_1}(a\!\circ\!\overline{\Gamma_2})
  = G_{\Gamma_1}(a)\!\circ\! G_{\Gamma_1}(\overline{\Gamma_2})
  = a\!\circ\!\overline{\Gamma_1}\!\circ\! G_{\Gamma_1}(\overline{\Gamma_2}).
\end{equation}
To compute $G_{\Gamma_1}(\overline{\Gamma_2})$ we need to write
$\overline{\Gamma_2}$ in terms of the generators of $\PCDc$. This
we do by replacing every $t^{12}$ appearing in $\Gamma_2$ by
$\overline{t^{12}}=a^{-1}\circ d_3H\circ a$, every $t^{13}$ by
$\overline{t^{13}}=d_0X\circ a^{-1}\circ d_3H\circ a \circ
d_0X$, and every $t^{23}$ by $\overline{t^{23}}=d_0H$. By the
definition of the action of $G_{\Gamma_1}$ on the generators
of $\PCDc$, we find that it maps $\overline{t^{12}}$ to
$\overline{\Gamma_1^{-1}t^{12}\Gamma_1}$, $\overline{t^{13}}$
to $\overline{(\Gamma_1^{-1})^{132}t^{13}\Gamma_1^{132}}$ and
$\overline{t^{23}}$ to $\overline{t^{23}}$. Combining this
and~\eqref{GammaComposition} we get~\eqref{GammaProduct}.
\end{pf}

Similar reasoning leads to the following:
\begin{proposition} The action of $\GRTc$ on $\ASSc$, written in terms of
$\Gamma$'s and $\Phi$'s, is given by
\[ \Gamma(\Phi)=\Gamma\cdot\left(\left.\Phi\right|_{
    t^{12}\to\Gamma_1^{-1}t^{12}\Gamma_1,\
    t^{13}\to(\Gamma_1^{-1})^{132}t^{13}\Gamma_1^{132},\
    t^{23}\to t^{23}
  }\right),
\]
with products and inverses taken in $\Apbc_3$. \qed
\end{proposition}

The group $\GTc$ admits a similar description. Any element of $\GTc$
maps $a$ to a limit of formal sums of parenthesized braids whose
skeleton is $a$. Such a limit is of the form $a\circ\Sigma$, where
$\Sigma$ is a self-morphism whose skeleton is the identity of the
object $(\bullet(\bullet\bullet))$ of $\PBc$, regarded as an element of
$\widehat{PB}_3$. Let $\sigma_1$ and $\sigma_2$ be the standard generators
$\setlength{\unitlength}{0.33\standardunitlength}
	\begin{array}{c}  \hspace{-1.7mm}
        	\raisebox{-2pt}{\input draws/sigma1.tex }
        	\hspace{-1.9mm}
	\end{array}
$ and $\setlength{\unitlength}{0.33\standardunitlength}
	\begin{array}{c}  \hspace{-1.7mm}
        	\raisebox{-2pt}{\input draws/sigma2.tex }
        	\hspace{-1.9mm}
	\end{array}
$ of the
(non-pure) braid group $B_3$ on 3 strands. Every $\Sigma\in\widehat{PB}_3$
is a limit of formal sums of combinations of $\sigma_{1,2}$.

\begin{proposition}
\begin{enumerate}
\item $\GTc$ can be identified as the group of all group-like
  non-degenerate $\Sigma\in\widehat{PB}_3$ satisfying:
  \begin{itemize}
  \item The pentagon for pure braids, in $\widehat{PB}_4$:
    \[ d_4\Sigma\cdot d_2\Sigma\cdot d_0\Sigma
       = d_1\Sigma\cdot d_3\Sigma
    \]
    (with the obvious interpretation for the $d_i$'s).
  \item The hexagons for pure braids, in $\hat{B}_3$, the unipotent
    completion of $B_3$:
    \[ \sigma_2\sigma_1=\Sigma\cdot\sigma_2\cdot\Sigma^{-1}\cdot
         \sigma_1\cdot\Sigma.
    \]
  \end{itemize}
\item The group law is given by
  \[ \Sigma_1\times\Sigma_2 = \Sigma_1\cdot\left(\left.\Sigma_2\right|_{
      \sigma_1\to\Sigma^{-1}\sigma_1\Sigma,\ 
      \sigma_2\to\sigma_2
    }\right),
  \]
  with products and inverses taken in $\hat{B}_3$.
\item The action on $\ASSc$ is given by
  \[ (\Phi,\Sigma)\mapsto\Phi^\Sigma
    = \Phi\cdot\left(\left.\Sigma\right|_{
      \sigma_1\to\Phi^{-1}e^{t^{12}/2}X_1\Phi,\ 
      \sigma_2\to e^{t^{23}/2}X_2
    }\right).
  \]
  This formula makes sense in $\Apbc_3\rtimes S_3$, with $X_1=(12)$ and
  $X_2=(23)$ the standard generators of the permutation group $S_3$ which
  acts on $\Apbc_3$ as in Definition~\ref{PermutationAction}. Implicitly
  we claim that this formula is well defined and valued in
  $\Apbc_3\subset\Apbc_3\rtimes S_3$. \qed
\end{enumerate}
\end{proposition}

%% file: refs.tex
\ifx\undefined\bysame
        \newcommand{\bysame}{\leavevmode\hbox to3em{\hrulefill}\,}
\fi

%% file: GT1.bbl
\begin{thebibliography}{B-N99}

\bibitem[B-N1]{Bar-Natan:Vassiliev} D.~Bar-Natan,
  {\em On the Vassiliev knot invariants},
  Topology {\bf 34} 423--472 (1995).

\bibitem[B-N2]{Bar-Natan:Homotopy} \bysame,
    {\em Vassiliev homotopy string link invariants},
  Jour.\ of Knot Theory and its Ramifications {\bf 4} (1995) 13--32.

\bibitem[B-N3]{Bar-Natan:NAT} \bysame,
    {\em Non-associative tangles},
  in {\em Geometric topology} (proceedings of the Georgia international
  topology conference), (W.~H.~Kazez, ed.), 139--183, Amer.\ Math.\
  Soc.\ and International Press, Providence, 1997.

\bibitem[B-N4]{Bar-Natan:Braids} \bysame,
  {\em Vassiliev and quantum invariants of braids},
  in Proc.\ of Symp.\ in Appl.\ Math.\ {\bf 51} (1996) 129--144 {\em The
  interface of knots and physics}, (L.~H.~Kauffman, ed.), Amer.\ Math.\
  Soc., Providence.

\bibitem[B-N5]{Bar-Natan:VasBib} \bysame,
        {\em Bibliography of Vassiliev Invariants},
        {\tt http://www.ma.huji.ac.il/$\sim$drorbn}.

\bibitem[Bi]{Birman:Bulletin} J.~S.~Birman,
  {\em New points of view in knot theory},
  Bull.\ Amer.\ Math.\ Soc.\ {\bf 28} (1993) 253--287.

\bibitem[BL]{BirmanLin:Vassiliev} \bysame and X-S.~Lin,
  {\em Knot polynomials and Vassiliev's invariants},
  Inv.\ Math.\ {\bf 111} (1993) 225--270.

\bibitem[Ca]{Cartier:Construction} P.~Cartier,
  {\em Construction combinatoire des invariants de Vassiliev-Kontsevich
    des n\oe uds},
  C.~R.~Acad. Sci. Paris {\bf 316} S\'{e}rie I (1993) 1205--1210.

\bibitem[Dr1]{Drinfeld:QuasiHopf} V.~G. Drinfel'd,
  {\em Quasi-Hopf algebras},
  Leningrad Math. J. {\bf 1} (1990) 1419--1457.

\bibitem[Dr2]{Drinfeld:GalQQ} \bysame,
  {\em On quasitriangular Quasi-Hopf algebras and a group closely
   connected with $\text{Gal}(\bar{{\Bbb Q}}/{\Bbb Q})$},
  Leningrad Math. J. {\bf 2} (1991) 829--860.

\bibitem[EK1]{EtingofKazhdan:Quantization} P.~Etingof and D.~Kazhdan,
    {\em Quantization of Lie bialgebras, I},
  Selecta Math.\ {\bf 2-1} (1996) 1--41. See also q-alg/9506005.

\bibitem[EK2]{EtingofKazhdan:Poisson} \bysame\ and \bysame,
    {\em Quantization of Poisson algebraic groups and Poisson homogeneous
      spaces},
  q-alg/9510020 preprint.

\bibitem[Go1]{Goussarov:New} M.~Goussarov,
  {\em A new form of the Conway-Jones polynomial of oriented links},
  in {\em Topology of manifolds and varieties} (O.~Viro, editor),
  Amer.\ Math.\ Soc., Providence 1994, 167--172.

\bibitem[Go2]{Goussarov:nEquivalence} \bysame,
  {\em On $n$-equivalence of knots and invariants of finite degree},
  in {\em Topology of manifolds and varieties} (O.~Viro, editor),
  Amer.\ Math.\ Soc., Providence 1994, 173--192.

\bibitem[Hu]{Hutchings:SingularBraids} M.~Hutchings,
  {\em Integration of singular braid invariants and graph cohomology},
  Harvard University preprint, April 1995.

\bibitem[Ka]{Kassel:Book} C.~Kassel,
  {\em Quantum groups},
  Springer-Verlag GTM {\bf 155}, New York 1994.

\bibitem[KT]{KasselTuraev:Graphs} \bysame\ and V.~Turaev,
  {\em Chord diagram invariants of tangles and graphs},
  University of Strasbourg preprint, January 1995.

\bibitem[Koh1]{Kohno:MonRep} T. Kohno,
  {\em Monodromy representations of braid groups and Yang-Baxter equations},
  Ann. Inst. Fourier {\bf 37} (1987) 139--160.

\bibitem[Koh2]{Kohno:deRham} \bysame,
  {\em Vassiliev invariants and de-Rham complex on the space of
   knots},
  Contemp.\ Math.\ {\bf 179} (1994) 123--138.

\bibitem[Kon]{Kontsevich:Vassiliev} M.~Kontsevich,
  {\em Vassiliev's knot invariants},
  Adv.\ in Sov.\ Math., {\bf 16(2)} (1993) 137--150.

\bibitem[Le]{Le:UniversalIHS} T.~Q.~T.~Le,
    {\em An invariant of integral homology 3-spheres which is universal
      for all finite type invariants},
  in {\em Solitons, geometry and topology: on the crossroad},
  (V.~Buchstaber and S.~Novikov, eds.) AMS Translations Series 2,
  Providence. See also q-alg/9601002.

\bibitem[LM1]{LeMurakami:Universal} \bysame\ and J.~Murakami,
    {\em The universal Vassiliev-Kontsevich invariant for framed
      oriented links},
  Compositio Math.\ {\bf 102} (1996), 42--64. See also hep-th/9401016.

\bibitem[LM2]{LeMurakami:Parallel} \bysame\ and \bysame,
  {\em Parallel version of the universal Vassiliev-Kontsevich invariant},
  to appear in J.\ Pure and Appl.\ Algebra.

\bibitem[LMO]{LeMurakamiOhtsuki:Universal} \bysame, \bysame\ and T.~Ohtsuki,
    {\em On a universal quantum invariant of 3-manifolds},
  to appear in Topology. See also q-alg/9512002.

\bibitem[LS1]{LochakSchneps:Automorphisms} P.~Lochak and L.~Schneps,
  {\em The Grothendieck-Teichmuller group and automorphisms of braid
    groups}, 
  in \cite{Schneps:DessinsdEnfants}, 323--358.

\bibitem[LS2]{LochakSchneps:Tower} \bysame and \bysame,
  {\em The $k$-pro-unipotent braid tower and the Grothendieck-Teichmuller
    group},
  preprint.

\bibitem[Ma]{MacLane:Categories} S.~Mac~Lane,
  {\em Categories for the working mathematician},
  Springer-Verlag GTM {\bf 5}, New-York 1971.

\bibitem[Pi]{Piunikhin:Combinatorial} S.~Piunikhin,
  {\em Combinatorial expression for universal Vassiliev link invariant},
  Commun.\ Math.\ Phys. {\bf 168-1} 1--22.

\bibitem[Sc]{Schneps:DessinsdEnfants} L.~Schneps (ed.),
  {\em The Grothendieck theory of dessins d'enfants},
  London Math.\ Soc.\ Lec.\ Notes Ser.\ {\bf 200} (1994), Cambridge Univ.\
  Press, Cambridge UK.

\bibitem[SS]{ShniderSternberg:Book} S.~Shnider and S.~Sternberg,
  {\em Quantum groups --- from coalgebras to Drinfel'd algebras},
  International Press, Cambridge MA 1994.

\bibitem[Va1]{Vassiliev:CohKnot} V.~A. Vassiliev,
  {\em Cohomology of knot spaces}, Theory of Singularities and
   its Applications (Providence) (V.~I. Arnold, ed.), Amer. Math.  Soc.,
  Providence, 1990.

\bibitem[Va2]{Vassiliev:Book} \bysame,
  {\em Complements of discriminants of smooth maps: topology and
   applications},
  Trans. of Math. Mono. {\bf 98}, Amer. Math. Soc., Providence, 1992.

\bibitem[Wa]{Waterhouse:GroupSchemes} W.~C.~Waterhouse,
  {\em Introduction to affine group schemes},
  Springer-Verlag GTM {\bf 66}, New-York 1979.

\bibitem[Za]{Zagier:Values} D.~Zagier,
  {\em Values of zeta functions and their applications},
  Max-Planck-Institut Bonn preprint.

\end{thebibliography}
